\documentclass[a4paper, 12pt, leqno]{article}
\usepackage{amsmath}
\usepackage{amscd}
\usepackage{amssymb}
\usepackage{amsthm}
\newtheorem{lemma}{{\sc Lemma}}[section]

\newtheorem{proposition}[lemma]{{\sc Proposition}}
\newtheorem{theorem}[lemma]{{\sc Theorem}}
\newtheorem{remark}[lemma]{{\sc Remark}}

\numberwithin{equation}{section}

\def\Ga{{\mathfrak{a}}}
\def\Gb{{\mathfrak{b}}}
\def\Gc{{\mathfrak{c}}}
\def\Gg{{\mathfrak{g}}}
\def\Gh{{\mathfrak{h}}}
\def\Gk{{\mathfrak{k}}}

\def\Gm{{\mathfrak{m}}}
\def\Gn{{\mathfrak{n}}}

\def\Gs{{\mathfrak{s}}}

\def\BA{{\mathbb{A}}}
\def\BC{{\mathbb{C}}}
\def\BF{{\BQ(q)}}
\def\BK{{\mathbb{K}}}

\def\BQ{{\mathbb{Q}}}
\def\BZ{{\mathbb{Z}}}

\def\CR{{R}}

\def\Hom{\mathop{\rm Hom}\nolimits}
\def\id{\mathop{\rm id}\nolimits}
\def\Image{\mathop{\rm Im}\nolimits}

\def\res{\mathop{\rm res}\nolimits}

\title{Poisson Hopf algebras associated to quantized enveloping algebras}
\date{}
\author{Toshiyuki TANISAKI
\cr\cr
{\footnotesize Department of Mathematics, 
Osaka City University,} \cr
{\footnotesize 3-3-138, Sugimoto, Sumiyoshi-ku,}\cr
{\footnotesize Osaka, 558-8585 Japan}
}
\begin{document}
\maketitle
\renewcommand{\thefootnote}{}
%\footnotetext
%{2000 {\it Mathematics Subject Classification:}
%Primary 20G42; Secondary 16S32, 17B37.
%}
%\footnotetext
%{{\it Key words and Phrases:}
%Poisson Hopf algebra, quantized enveloping algebra.
%}
\begin{abstract}
We study certain Poisson structures related to quantized enveloping algebras.
In particular, we give a description of the Poisson structure of a certain manifold associated to the ring of differential operators.
\end{abstract}
%\tableofcontents
\setcounter{section}{-1}
\section{Introduction}
\label{sec:Intoro}
Let $U$ be the quantized enveloping algebra corresponding to a finite dimensional complex simple Lie algebra $\Gg$.
It is a Hopf algebra over $\BQ(q,q^{-1})$.
In \cite{DP} De Concini-Procesi introduced a certain $\BQ[q,q^{-1}]$-form $U_{\BQ[q,q^{-1}]}$, and studied its properties (see also \cite{DeConcini}).
For $z\in\BC^\times$ we denote by 
$U_z$
the specialization of $U_{\BQ[q,q^{-1}]}$ at $q=z$.
Let $\ell$ be a positive odd integer, and let $\zeta\in\BC^\times$ be a primitive $\ell$-th root of $1$ (we assume that $\ell$ is prime to 3 for type $G_2$).
Then $U_\zeta$ is canonically isomorphic to the specialization of the De Concini-Kac form at $q=\zeta$.
Denote by $Z_\zeta$ the central Hopf subalgebra of $U_\zeta$ generated by the $\ell$-th powers of typical generators. 
It is called the Frobenius center of $U_\zeta$.
Then one of the main results of \cite{DP} (and \cite{DeConcini}) is the following isomorphisms of Hopf algebras:
\begin{equation}
\label{eq:0}
Z_\zeta\cong U_1\cong \BC[M],
\end{equation}
where $M$ is a certain algebraic group associated to $\Gg$.
The three Hopf algebras appearing in \eqref{eq:0} are endowed with natural Poisson Hopf algebra structures, and the isomorphisms in \eqref{eq:0} is in fact that of Poisson Hopf algebras.
In \cite{DP} De Concini-Procesi constructed an isomorphism $U_1\cong \BC[M]$ by giving a correspondence between generators and verifying the necessary relations among generators by a lengthy calculation.
Later a more natural construction of the isomorphism in terms of the Drinfleld paring was found by Gavarini \cite{Gav}.

In this note we present a slightly different proof of \eqref{eq:0} which is still based on Gavarini's construction of the isomorphism.
In fact in \cite{Gav} the statement about the Poisson structure is deduced from its dual statement, but our argument is more direct (we do not claim that our proof is simpler than the one in \cite{Gav}).
We also give a description of the Poisson algebra associated to the ring of differential operators.

In this paper we shall use the following notation for a Hopf algebra $H$ over a field $\BK$.
The comultiplication, the counit, and the antipode of $H$ are denoted by
\begin{align}
&\Delta_H:H\to H\otimes_\BK H,\\
&\varepsilon_H:H\to\BK,\\
&S_H:H\to H
\end{align}
respectively.
The subscript $H$ will often be omitted.
For $n\in\BZ_{>0}$ we denote by
\[
\Delta_n:H\to H^{\otimes n+1}
\]
the algebra homomorphism given by 
\[
\Delta_1=\Delta,\qquad
\Delta_n=(\Delta\otimes \id_{H^{\otimes n-1}})\circ\Delta_{n-1},
\]
and write
\[
\Delta(h)=\sum_{(h)}h_{(0)}\otimes h_{(1)},
\qquad
\Delta_n(h)=\sum_{(h)_n}h_{(0)}\otimes\cdots\otimes h_{(n)}
\quad(n\geqq2).
\]

I would like to express my appreciation to F. Gavarini for 
pointing out the reference \cite{Gav} and informing me of 
several valuable comments on the first version of this manuscript.

\section{Lie algebras}
Let $\Gg$ be a finite-dimensional simple Lie algebra over $\BC$ and let $\Gh$ be its Cartan subalgebra.
We denote by $\Delta\subset\Gh^*$, $Q\subset\Gh^*$ and $W\subset GL(\Gh^*)$ the set of roots, the root lattice $\sum_{\alpha\in\Delta}\BZ\alpha$, and the Weyl group respectively.
We fix a set of simple roots $\{\alpha_i\}_{i\in I}$, and denote the corresponding set of positive roots and simple reflections by $\Delta^+\subset\Gh^*$ and  $\{s_i\}_{i\in I}\subset W$ respectively.
Set 
\[
Q^+
=\sum_{\alpha\in\Delta^+}\BZ_{\geqq0}\alpha
=\bigoplus_{i\in I}\BZ_{\geqq0}\alpha_i
\subset \Gh^*.
\]
We denote the longest element of $W$ by $w_0$.
Let
\begin{equation}
\label{eq:SBMh*}
(\,,\,):\Gh^*\times\Gh^*\to\BC
\end{equation}
be the $W$-invariant symmetric bilinear form 
satisfying $(\beta,\beta)/2=1$ for short roots $\beta\in\Delta$.
We define subalgebras $\Gn^+, \Gn^-$ of $\Gg$ by
\[
\Gn^\pm
=\bigoplus_{\beta\in\Delta^+}\Gg_{\pm\beta},
\]
where
\[
\Gg_{\pm\beta}
=\{x\in\Gg\mid [h,x]=\pm\beta(h)x\,\,(h\in \Gh)\}.
\]
Then we have
\[
\Gg=\Gn^+\oplus\Gh\oplus\Gn^-
\]
For each $i\in I$ we take ${e}_i\in\Gg_{\alpha_i}, {f}_i\in\Gg_{-\alpha_i}, h_i\in\Gh$ such that $[{e}_i,{f}_i]=h_i$ and $\alpha_i(h_i)=2$.

We define a subalgebra $\Gm$ of $\Gg\oplus\Gg$ by 
\[
\Gm=\{(h+x,-h+y)\mid h\in\Gh, x\in\Gn^+, y\in\Gn^-\}.
\]
Set 
\[
\Gm^0=\{(h,-h)\mid h\in\Gh\},\quad
\Gm^+=\{(x,0)\mid x\in\Gn^+\},\quad
\Gm^-=\{(0,y)\mid y\in \Gn^-\}.
\]
They are subalgebras of $\Gm$, and we have $\Gm=\Gm^+\oplus\Gm^0\oplus\Gm^-$.
Moreover, we have isomorphisms
\begin{align}
\label{ident1}
&\Gh\simeq\Gm^0\qquad(h\leftrightarrow(h,-h)),\\
\label{ident2}
&\Gn^+\simeq\Gm^+\qquad(x\leftrightarrow(x,0)),\\
\label{ident3}
&\Gn^-\simeq\Gm^-\qquad(y\leftrightarrow(0,y))
\end{align}
of Lie algebras.
We denote by 
\begin{equation}
\iota:\Gh\to\Gm^0
\end{equation}
the isomorphism \eqref{ident1}.
For $i\in I$ set
\[
x_i=(e_i,0)\in\Gm^+,\quad y_i=(0, f_i)\in\Gm^-,\quad
t_i=(
h_i,-h_i)\in\Gm^0.
\]
Let $G$ be the adjoint group of $\Gg$.
We denote by $M$ the connected closed subgroup of $G\times G$ with Lie algebra $\Gm$.
Let $M^0$, $M^\pm$ be the connected closed subgroups of $M$ with Lie algebras $\Gm^0$, $\Gm^\pm$ respectively.

\section{Quantized enveloping algebras}
For $n\in\BZ_{\geqq0}$ we set
\[
[n]_t=\frac{t^n-t^{-n}}{t-t^{-1}}\in\BZ[t,t^{-1}],
\qquad
[n]_t!=[n]_t[n-1]_t\cdots[2]_t[1]_t\in\BZ[t,t^{-1}].
\]
For $n\in\BZ$ and $m\in\BZ_{\geqq0}$ we set
\[
\begin{bmatrix}n\\m\end{bmatrix}_t=
[n]_t[n-1]_t\cdots[n-m+1]_t/[m]_t!.
\]

The quantized enveloping algebra $U=U_q(\Gg)$ of $\Gg$ is an associative algebra over  $\BF=\BQ(q)$ with identity element $1$ generated by the elements $K_\lambda\,(\lambda\in Q),\,E_i, F_i\,(i\in I)$ satisfying the following defining relations:
\begin{align}
&K_0=1,\quad 
K_\lambda K_\mu=K_{\lambda+\mu}
\qquad(\lambda,\mu\in Q),
\label{eq:def1}\\
&K_\lambda E_iK_\lambda^{-1}=q^{(\lambda,\alpha_i)}E_i,\qquad(\lambda\in Q, i\in I),
\label{eq:def2a}\\
&K_\lambda F_iK_\lambda^{-1}=q^{-(\lambda,\alpha_i)}F_i
\qquad(\lambda\in Q, i\in I),
\label{eq:def2b}\\
&E_iF_j-F_jE_i=\delta_{ij}\frac{K_i-K_i^{-1}}{q_i-q_i^{-1}}
\qquad(i, j\in I),
\label{eq:def3}\\
&\sum_{n=0}^{1-a_{ij}}(-1)^nE_i^{(1-a_{ij}-n)}E_jE_i^{(n)}=0
\qquad(i,j\in I,\,i\ne j),
\label{eq:def4}\\
&\sum_{n=0}^{1-a_{ij}}(-1)^nF_i^{(1-a_{ij}-n)}F_jF_i^{(n)}=0
\qquad(i,j\in I,\,i\ne j),
\label{eq:def5}
\end{align}
where $q_i=q^{(\alpha_i,\alpha_i)/2}, K_i=K_{\alpha_i}, a_{ij}=2(\alpha_i,\alpha_j)/(\alpha_i,\alpha_i)$ for $i, j\in I$, and
\[
E_i^{(n)}=
E_i^n/[n]_{q_i}!,
\qquad
F_i^{(n)}=
F_i^n/[n]_{q_i}!
\]
for $i\in I$ and $n\in\BZ_{\geqq0}$.
Algebra homomorphisms $\Delta:U\to U\otimes U, \varepsilon:U\to\BF$ and an algebra anti-automorphism $S:U\to U$ are defined by:
\begin{align}
&\Delta(K_\lambda)=K_\lambda\otimes K_\lambda,\\
&\Delta(E_i)=E_i\otimes 1+K_i\otimes E_i,\quad
\Delta(F_i)=F_i\otimes K_i^{-1}+1\otimes F_i,
\nonumber\\
&\varepsilon(K_\lambda)=1,\quad
\varepsilon(E_i)=\varepsilon(F_i)=0,\\
&S(K_\lambda)=K_\lambda^{-1},\quad
S(E_i)=-K_i^{-1}E_i, \quad S(F_i)=-F_iK_i,
\end{align}
and $U$ is endowed with a Hopf algebra structure with the comultiplication $\Delta$, the counit $\varepsilon$ and the antipode $S$.

We define subalgebras $U^0, U^{\geqq0}, U^{\leqq0}, U^{+}, U^{-}$ of $U$ by 
\begin{align}
U^0&=\langle K_\lambda\mid \lambda\in Q\rangle,\\
U^{\geqq0}&=\langle K_\lambda, E_i\mid \lambda\in Q, i\in I\rangle,\\
U^{\leqq0}&=\langle K_\lambda, F_i\mid \lambda\in Q, i\in I\rangle,\\
U^{+}&=\langle E_i\mid i\in I\rangle,\\
U^{-}&=\langle F_i\mid i\in I\rangle.
\end{align}

Note that $U^0, U^{\geqq0}, U^{\leqq0}$ are Hopf subalgebras of $U$, while $U^+$ and $U^-$ are not Hopf subalgebras.

The following result is standard.
\begin{proposition}
\label{prop:Str-of-U}
\begin{itemize}
\item[\rm(i)]
$\{K_\lambda\mid \lambda\in Q\}$ is a $\BF$-basis of $U^0$.
\item[\rm(ii)]
$U^+$ $($resp.\ $U^-$$)$ is isomorphic to the $\BF$-algebra generated by 
$\{E_i\mid i\in I\}$ $($resp.\ $\{F_i\mid i\in I\}$$)$ with defining relation 
\eqref{eq:def4} $($resp.\ \eqref{eq:def5}$)$.
\item[\rm(iii)]
$U^{\geqq0}$ $($resp.\ $U^{\leqq0}$$)$ is isomorphic to the $\BF$-algebra generated by $\{E_i, K_\lambda\mid i\in I, \lambda\in Q\}$ 
$($resp.\ $\{F_i, K_\lambda\mid i\in I, \lambda\in Q\}$$)$ 
with defining relations \eqref{eq:def1}, \eqref{eq:def2a}, \eqref{eq:def4} $($resp.\ \eqref{eq:def1}, \eqref{eq:def2b}, \eqref{eq:def5}$)$.
\item[\rm(iv)]
The linear maps 
\begin{gather*}
U^-\otimes U^0\otimes U^+\to U \gets U^+\otimes U^0\otimes U^-,\\
U^+\otimes U^0\to U^{\geqq0} \gets U^0\otimes U^+,\qquad
U^-\otimes U^0\to U^{\leqq0} \gets U^0\otimes U^-
\end{gather*}
induced by the multiplication are all isomorphisms.
\end{itemize}
\end{proposition}

For $\gamma\in Q$ we set
\[
U^\pm_{\gamma}
=\{x\in U^\pm\mid
K_\lambda xK_\lambda^{-1}=q^{(\lambda,\gamma)}x\,\,(\lambda\in Q)\}.
\]
We have $U^\pm_{\pm\gamma}=\{0\}$ unless $\gamma\in Q^+$, and 
\[
U^\pm=\bigoplus_{\gamma\in Q^+}U^\pm_{\pm\gamma}.
\]

For $i\in I$ we can define an algebra automorphism $T_i$ of $U$ by
\begin{align*}
&T_i(K_\mu)=K_{s_i\mu}\qquad(\mu\in Q),\\
&T_i(E_j)=
\begin{cases}
\sum_{k=0}^{-a_{ij}}(-1)^kq_i^{-k}E_i^{(-a_{ij}-k)}E_jE_i^{(k)}\qquad
&(j\in I,\,\,j\ne i),\\
-F_iK_i\qquad
&(j=i),
\end{cases}\\
&T_i(F_j)=
\begin{cases}
\sum_{k=0}^{-a_{ij}}(-1)^kq_i^{k}F_i^{(k)}F_jF_i^{(-a_{ij}-k)}\qquad
&(j\in I,\,\,j\ne i),\\
-K_i^{-1}E_i\qquad
&(j=i).
\end{cases}
\end{align*}
For $w\in W$ we define an algebra automorphism $T_w$ of $U$ by $T_w=T_{i_1}\cdots T_{i_n}$ where $w=s_{i_1}\cdots s_{i_n}$ is a reduced expression.
The automorphism $T_w$ does not depend on the choice of a reduced expression (see Lusztig \cite{Lbook}).

We fix a reduced expression 
\[
w_0=s_{i_1}\cdots s_{i_N}
\]
of $w_0$, and set
\[
\beta_k=s_{i_1}\cdots s_{i_{k-1}}(\alpha_{i_k})\qquad
(1\leqq k\leqq N).
\]
Then we have $\Delta^+=\{\beta_k\mid1\leqq k\leqq N\}$.
For $1\leqq k\leqq N$ set 
\begin{equation}
\label{eq:root vector}
E_{\beta_k}=T_{i_1}\cdots T_{i_{k-1}}(E_{i_k}),\quad
F_{\beta_k}=T_{i_1}\cdots T_{i_{k-1}}(F_{i_k}).
\end{equation}
Then $\{E_{\beta_{N}}^{m_N}\cdots E_{\beta_{1}}^{m_1}\mid
m_1,\dots, m_N\geqq0\}$ (resp. $\{F_{\beta_{N}}^{m_N}\cdots F_{\beta_{1}}^{m_1}\mid
m_1,\dots, m_N\geqq0\}$)
is a $\BF$-basis of $U^+$ (resp. $U^-$), called the PBW-basis (see Lusztig \cite{L2}).
We have $E_{\alpha_i}=E_i$ and $F_{\alpha_i}=F_i$ for any $i\in I$.
For $1\leqq k\leqq N,\, m\geqq0$ we also set 
\begin{equation}
E^{(m)}_{\beta_k}=E^{m}_{\beta_k}/[m]_{q_{\beta_k}}!,\quad
F^{(m)}_{\beta_k}=F^{m}_{\beta_k}/[m]_{q_{\beta_k}}!,
\end{equation}
where $q_\beta=q^{(\beta,\beta)/2}$ for $\beta\in\Delta^+$.

There exists a unique bilinear form
\begin{equation}
\label{eq:Drinfeld-paring}
\tau:U^{\geqq0}\times U^{\leqq0}\to\BF,
\end{equation}
called the Drinfeld paring, which is characterized by
\begin{align}
&\tau(x,y_1y_2)=(\tau\otimes\tau)(\Delta(x),y_1\otimes y_2)
&(x\in U^{\geqq0},\,y_1,y_2\in U^{\leqq0}),
\label{eq:DP1}\\
&\tau(x_1x_2,y)=(\tau\otimes\tau)(x_2\otimes x_1,\Delta(y))
&(x_1, x_2\in U^{\geqq0},\,y\in U^{\leqq0}),
\label{eq:DP2}\\
&\tau(K_\lambda,K_\mu)=q^{-(\lambda,\mu)}
&(\lambda,\mu\in Q),
\label{eq:DP3}\\
&\tau(K_\lambda, F_i)=\tau(E_i,K_\lambda)=0
&(\lambda\in Q,\,i\in I),
\label{eq:DP4}\\
&\tau(E_i,F_j)=\delta_{ij}/(q_i^{-1}-q_i)
&(i,j\in I).
\label{eq:DP5}
\end{align}
It satisfies the following (see Tanisaki \cite{T}, Lusztig \cite{Lbook}).
\begin{lemma}
\label{lem:Drinfeld paring}
\begin{itemize}
\item[\rm(i)]
$\tau(S(x),S(y))=\tau(x,y)$ for $x\in U^{\geqq0}, y\in U^{\leqq0}$.
\item[\rm(ii)]
For $x\in U^{\geqq0}, y\in U^{\leqq0}$ we have
\begin{align*}
yx=\sum_{(x)_2,(y)_2}
\tau(x_{(0)},S(y_{(0)}))\tau(x_{(2)},y_{(2)})x_{(1)}y_{(1)},\\
xy=\sum_{(x)_2,(y)_2}
\tau(x_{(0)},y_{(0)})\tau(x_{(2)},S(y_{(2)}))y_{(1)}x_{(1)}.
\end{align*}
\item[\rm(iii)]
$\tau(xK_\lambda, yK_\mu)=q^{-(\lambda,\mu)}\tau(x,y)$ for $\lambda, \mu\in Q, x\in U^+, y\in U^-$.
\item[\rm(iv)]
$\tau(U^+_\beta, U^-_{-\gamma})=\{0\}$ for $\beta, \gamma\in Q^+$ with $\beta\ne\gamma$.
\item[\rm(v)]
For any $\beta\in Q^+$ the restriction of \eqref{eq:Drinfeld-paring} to $U^+_\beta\times U^-_{-\beta}$ is non-degenerate.
\end{itemize}
\end{lemma}
We have the following explicit computation of \eqref{eq:Drinfeld-paring} in terms of PBW-bases (\cite{KR}, \cite{KT}, \cite{LS}).
\begin{proposition}
\label{prop:DPPBW}
We have
\begin{align*}
&\tau(E_{\beta_{N}}^{m_N}\cdots E_{\beta_{1}}^{m_1},
F_{\beta_{N}}^{n_N}\cdots F_{\beta_{1}}^{n_1})\\
=&
\prod_{k=1}^N
\delta_{m_k,n_k}
(-1)^{m_k}
[m_k]_{q_{\beta_k}}!
q_{\beta_k}^{m_k(m_k-1)/2}
(q_{\beta_k}-q_{\beta_k}^{-1})^{-m_k}.
\end{align*}
\end{proposition}

The quantized enveloping algebra $V=U_q(\Gm)$ of $\Gm$ is an associative algebra over  $\BF$ with identity element $1$ generated by the elements $Z_\lambda\,(\lambda\in Q),\,X_i, Y_i\,(i\in I)$ satisfying the following defining relations:
\begin{align}
&Z_0=1,\quad 
Z_\lambda Z_\mu=Z_{\lambda+\mu}
\qquad(\lambda,\mu\in Q),
\label{eq:def1V}\\
&Z_\lambda X_iZ_\lambda^{-1}=q^{(\lambda,\alpha_i)}X_i,\qquad(\lambda\in Q, i\in I),
\label{eq:def2aV}\\
&Z_\lambda Y_iZ_\lambda^{-1}=q^{(\lambda,\alpha_i)}Y_i
\qquad(\lambda\in Q, i\in I),
\label{eq:def2bV}\\
&X_iY_j-Y_jX_i=0
\qquad(i, j\in I),
\label{eq:def3V}\\
&\sum_{n=0}^{1-a_{ij}}(-1)^nX_i^{(1-a_{ij}-n)}X_jX_i^{(n)}=0
\qquad(i,j\in I,\,i\ne j),
\label{eq:def4V}\\
&\sum_{n=0}^{1-a_{ij}}(-1)^nY_i^{(1-a_{ij}-n)}Y_jY_i^{(n)}=0
\qquad(i,j\in I,\,i\ne j),
\label{eq:def5V}
\end{align}
where 
\[
X_i^{(n)}=
X_i^n/[n]_{q_i}!,
\qquad
Y_i^{(n)}=
Y_i^n/[n]_{q_i}!.
\]
$V$ is endowed with a structure of Hopf algebra by 
\begin{align}
&\Delta(Z_\lambda)=Z_\lambda\otimes Z_\lambda,\\
&\Delta(X_i)=X_i\otimes 1+Z_i\otimes X_i,\quad
\Delta(Y_i)=Y_i\otimes Z_i+1\otimes Y_i,
\nonumber\\
&\varepsilon(Z_\lambda)=1,\quad
\varepsilon(X_i)=\varepsilon(Y_i)=0,\\
&S(Z_\lambda)=Z_\lambda^{-1},\quad
S(X_i)=-Z_i^{-1}X_i, \quad S(Y_i)=-Y_iZ_i^{-1},
\end{align}
where $Z_i=Z_{\alpha_i}$ for $i\in I$.

We define subalgebras $V^0, V^{\geqq0}, V^{\leqq0}, V^{+}, V^{-}$ of $V$ by 
\begin{align}
V^0&=\langle Z_\lambda\mid \lambda\in Q\rangle,\\
V^{\geqq0}&=\langle Z_\lambda, X_i\mid \lambda\in Q, i\in I\rangle,\\
V^{\leqq0}&=\langle Z_\lambda, Y_i\mid \lambda\in Q, i\in I\rangle,\\
V^{+}&=\langle X_i\mid i\in I\rangle,\\
V^{-}&=\langle Y_i\mid i\in I\rangle.
\end{align}
Then $V^0$, $V^{\geqq0}$, $V^{\leqq0}$ are Hopf subalgebras of $V$.

Similarly to Proposition \ref{prop:Str-of-U} we have the following.
\begin{proposition}
\label{prop:Str-of-V}
\begin{itemize}
\item[\rm(i)]
$\{Z_\lambda\mid \lambda\in Q\}$ is a $\BF$-basis of $V^0$.
\item[\rm(ii)]
$V^+$ $($resp.\ $V^-$$)$ is isomorphic to the $\BF$-algebra generated by 
$\{X_i\mid i\in I\}$ $($resp.\ $\{Y_i\mid i\in I\}$$)$ with defining relation 
\eqref{eq:def4V} $($resp.\ \eqref{eq:def5V}$)$.
\item[\rm(iii)]
$V^{\geqq0}$ $($resp.\ $V^{\leqq0}$$)$ is isomorphic to the $\BF$-algebra generated by $\{X_i, Z_\lambda\mid i\in I, \lambda\in Q\}$ 
$($resp.\ $\{Y_i, Z_\lambda\mid i\in I, \lambda\in Q\}$$)$ 
with defining relations \eqref{eq:def1V}, \eqref{eq:def2aV}, \eqref{eq:def4V} $($resp.\ \eqref{eq:def1V}, \eqref{eq:def2bV}, \eqref{eq:def5V}$)$.
\item[\rm(iv)]
The linear maps 
\begin{gather*}
V^-\otimes V^0\otimes V^+\to V \gets V^+\otimes V^0\otimes V^-,\\
V^+\otimes V^0\to V^{\geqq0} \gets V^0\otimes V^+,\qquad
V^-\otimes V^0\to V^{\leqq0} \gets V^0\otimes V^-
\end{gather*}
induced by the multiplication are all isomorphisms.
\end{itemize}
\end{proposition}

For $\gamma\in Q$ we set
\[
V^\pm_{\gamma}
=\{x\in V^\pm\mid
Z_\lambda xZ_\lambda^{-1}=q^{\pm(\lambda,\gamma)}x\,\,(\lambda\in Q)\}.
\]
We have $V^\pm_{\pm\gamma}=\{0\}$ unless $\gamma\in Q^+$, and
\[
V^\pm=\bigoplus_{\gamma\in Q^+}V^\pm_{\pm\gamma}.
\]

By Proposition \ref{prop:Str-of-U} and Proposition \ref{prop:Str-of-V} we have isomorphisms
\begin{align*}
&\eta^{\leqq0}:V^{\leqq0}\to U^{\leqq0}
\qquad(Y_i\mapsto F_i,\,Z_\lambda\mapsto K_{-\lambda}),\\
&\eta^{\geqq0}:V^{\geqq0}\to U^{\geqq0}
\qquad(X_i\mapsto E_i,\,Z_\lambda\mapsto K_{\lambda})
\end{align*}
of  Hopf algebras.

We define a bilinear form
\begin{equation}
\sigma:U\times V\to\BF
\end{equation}
by
\begin{multline*}
\sigma(u_+u_0S(u_-),v_-v_+v_0)=
\tau(u_+,\eta^{\leqq0}(v_-))
\tau(u_0,\eta^{\leqq0}(v_0))
\tau(\eta^{\geqq0}(v_+),u_-)\\
\qquad(u_\pm\in U^\pm, u_0\in U^0, v_\pm\in V^\pm, v_0\in V^0).
\end{multline*}
Note that 
\begin{multline*}
\sigma(u_+u_{\leqq0},v_{\geqq0})=\tau((S\circ\eta^{\geqq0})(v_{\geqq0}),u_{\leqq0})\varepsilon (u_{+})\\
(v_{\geqq0}\in V^{\geqq0}, u_+\in U^+, u_{\leqq0}\in U^{\leqq0}),
\end{multline*}
\begin{multline*}
\sigma(u_{\geqq0}S(u_-),v_{\leqq0})
=\tau(u_{\geqq0},\eta^{\leqq0}(v_{\leqq0}))\varepsilon (u_{-})\\
(v_{\leqq0}\in V^{\leqq0}, u_-\in U^-, u_{\geqq0}\in U^{\geqq0}).
\end{multline*}
The following result is a consequence of  Gavarini \cite[Theorem 6.2]{Gav}.
\begin{proposition}
\label{prop:invariance}
We have
\[
\sigma(u,vv')=(\sigma\otimes\sigma)(\Delta(u),v\otimes v')
\qquad(u\in U, v, v'\in V).
\]
\end{proposition}

\section{$\BA$-forms}
We fix a subring $\BA$ of $\BF$ containing $\BQ[q,q^{-1}]$.
We denote by $U_{\BA}^L$ the Lusztig ${\BA}$-form of $U$, i.e., $U_{\BA}^L$ is the ${\BA}$-subalgebra of  $U$ generated by the elements
\[
E_i^{(m)},\, F_i^{(m)},\, K_\lambda\qquad
(i\in I,\, m\geqq0,\,\lambda\in Q).
\]
Set 
\begin{align*}
&U_{\BA}^{L,\pm}=U_{\BA}^L\cap U^\pm,\qquad
U_{\BA}^{L,0}=U_{\BA}^L\cap U^0,\\
&U_{\BA}^{L,\geqq0}=U_{\BA}^L\cap U^{L,\geqq0},\qquad
U_{\BA}^{L,\leqq0}=U_{\BA}^L\cap U^{L,\leqq0},
\\
&U_{{\BA},\pm\gamma}^{L,\pm}=U_{\BA}^L\cap U_{\pm\gamma}^\pm\quad(\gamma\in Q^+).
\end{align*}
Then $U_{\BA}^{L}, U_{\BA}^{L,0}, U_{\BA}^{L,\geqq0}, U_{\BA}^{L,\leqq0}$ are endowed with structures of Hopf algebras over ${\BA}$ via the Hopf algebra structure on $U$, and the multiplication of $U_{\BA}^L$ induces isomorphisms
\begin{align*}
&U_{\BA}^{L}\simeq
U_{\BA}^{L,-}\otimes
U_{\BA}^{L,0}\otimes
U_{\BA}^{L,+}
\simeq
U_{\BA}^{L,+}\otimes
U_{\BA}^{L,0}\otimes
U_{\BA}^{L,-},\\
&U_{\BA}^{L,\geqq0}\simeq
U_{\BA}^{L,0}\otimes
U_{\BA}^{L,+}
\simeq
U_{\BA}^{L,+}\otimes
U_{\BA}^{L,0},\\
&U_{\BA}^{L,\leqq0}\simeq
U_{\BA}^{L,0}\otimes
U_{\BA}^{L,-}
\simeq
U_{\BA}^{L,-}\otimes
U_{\BA}^{L,0}
\end{align*}
of ${\BA}$-modules.
Moreover, $U_{\BA}^{L,+}$, $U_{\BA}^{L,-}$, $U_{\BA}^{L,0}$ are free ${\BA}$-modules with bases
\begin{align*}
&\{E_{\beta_{N}}^{(m_N)}\cdots E_{\beta_{1}}^{(m_1)}\mid
m_1,\dots, m_N\geqq0\},\\
&\{F_{\beta_{N}}^{(m_N)}\cdots F_{\beta_{1}}^{(m_!)}\mid
m_1,\dots, m_N\geqq0\},\\
&
\left\{
\prod_{i\in I}
\left(
K_i^{\varepsilon_i}
\begin{bmatrix}{K_i}\\{n_i}
\end{bmatrix}
\right)
\mid
n_i\geqq0, \,\,
\varepsilon_i=0 \mbox{ or } 1
\right\}
\end{align*}
respectively, 
where
\[
\begin{bmatrix}{K_i}\\{m}\end{bmatrix}
=\prod_{s=0}^{m-1}
\frac{q_i^{-s}K_i-q_i^sK_i^{-1}}
{q_i^{s+1}-q_i^{-s-1}}
\qquad(m\geqq0).
\]

We denote by $V_{\BA}$ the ${\BA}$-subalgebra of $V$ generated by the elements
\[
X_i^{(m)},\, Y_i^{(m)},\, Z_i^{\pm1},\,
 \begin{bmatrix}Z_i\\m\end{bmatrix}\qquad
(i\in I,\, m\geqq0),
\]
where
\[
\begin{bmatrix}{Z_i}\\{m}\end{bmatrix}
=\prod_{s=0}^{m-1}
\frac{q_i^{-s}Z_i-q_i^sZ_i^{-1}}
{q_i^{s+1}-q_i^{-s-1}}
\qquad(m\geqq0).
\]
Set 
\begin{align*}
&V_{\BA}^{\pm}=V_{\BA}\cap V^\pm,\qquad
V_{\BA}^{0}=V_{\BA}\cap V^0,\\
&V_{\BA}^{\geqq0}=V_{\BA}\cap V^{\geqq0},\qquad
V_{\BA}^{\leqq0}=V_{\BA}\cap V^{\leqq0},
\\
&V_{{\BA},\pm\gamma}^{\pm}=V_{\BA}\cap V_{\pm\gamma}^\pm\quad(\gamma\in Q^+).
\end{align*}
Then we have
\[
V^+_{\BA}=\langle X_i^{(m)}\mid i\in I, m\geqq0\rangle,\qquad
V^-_{\BA}=\langle Y_i^{(m)}\mid i\in I, m\geqq0\rangle
\]
as ${\BA}$-algebras, and $V^0_{\BA}$ is a free ${\BA}$-module with basis consisting of the elements
\[
\prod_{i\in I}\left(
Z_i^{\varepsilon_i}
\begin{bmatrix}{Z_i}\\{m_i}\end{bmatrix}
\right)
\qquad(m_i\geqq0,\,\,\varepsilon_i=0 \mbox{ or } 1).
\]
Moreover, the multiplication of $V_{\BA}$ induces isomorphisms
\begin{align*}
&V_{\BA}\simeq
V_{\BA}^{-}\otimes
V_{\BA}^{0}\otimes
V_{\BA}^{+}
\simeq
V_{\BA}^{+}\otimes
V_{\BA}^{0}\otimes
V_{\BA}^{-},\\
&V_{\BA}^{\geqq0}\simeq
V_{\BA}^{0}\otimes
V_{\BA}^{+}
\simeq
V_{\BA}^{+}\otimes
V_{\BA}^{0},\\
&V_{\BA}^{\leqq0}\simeq
V_{\BA}^{0}\otimes
V_{\BA}^{-}
\simeq
V_{\BA}^{-}\otimes
V_{\BA}^{0}
\end{align*}
of ${\BA}$-modules.
Note that 
\begin{align*}
&\eta^{\geqq0}(V_{\BA}^{\geqq0})=U^{L,\geqq0}_{\BA},\qquad
\eta^{\geqq0}(V_{\BA}^{+})=U^{L,+}_{\BA},\qquad
\eta^{\geqq0}(V_{\BA}^{0})=U^{L,0}_{\BA},
\\
&\eta^{\leqq0}(V_{\BA}^{\leqq0})=U^{L,\leqq0}_{\BA},\qquad
\eta^{\leqq0}(V_{\BA}^{-})=U^{L,-}_{\BA},\qquad
\eta^{\leqq0}(V_{\BA}^{0})=U^{L,0}_{\BA}.
\end{align*}

We define root vectors $X_{\beta_k}\in V^+_{\beta_k}, Y_{\beta_k}\in V^-_{-\beta_k}\,\,(1\leqq k\leqq N)$ by
\begin{equation}
\eta^{\geqq0}(X_{\beta_k})=E_{\beta_k},\qquad
\eta^{\leqq0}(Y_{\beta_k})=F_{\beta_k},
\end{equation}
and set
\begin{equation}
\label{eq:root-victor-XY}
X_{\beta_k}^{(m)}=
\frac{X_{\beta_k}^m}{[m]_{q_{\beta_k}}!},\qquad
Y_{\beta_k}^{(m)}=
\frac{Y_{\beta_k}^m}{[m]_{q_{\beta_k}}!}
\qquad
(1\leqq k\leqq N,\,\,m\geqq0).
\end{equation}
Then we have free bases
$\{X_{\beta_{N}}^{(m_N)}\cdots X_{\beta_{1}}^{(m_1)}\mid
m_1,\dots, m_N\geqq0\}$ and $\{Y_{\beta_{N}}^{(m_N)}\cdots Y_{\beta_{1}}^{(m_1)}\mid
m_1,\dots, m_N\geqq0\}$ of $V_{\BA}^+$ and $V_{\BA}^-$ respectively.

Set
\begin{align}
&U_{\BA}=\{u\in U\mid\sigma(u,V_{\BA})\subset {\BA}\},\\
&
U_{\BA}^\pm=U^\pm\cap U_{\BA},\qquad
U_{\BA}^{0}=U^0\cap U_{\BA},\\
&U_{\BA}^{\geqq0}=U^{\geqq0}\cap U_{\BA},\qquad
U_{\BA}^{\leqq0}=U^{\leqq0}\cap U_{\BA},\\
&U_{{\BA},\pm\beta}^{\pm}
=U_{{\BA}}^{\pm}\cap U_{\pm\beta}^{\pm}
\qquad(\beta\in Q^+).
\end{align}
We can easily check that
\begin{align}
\label{eq:UA+}
&U^+_{\BA}=\{x\in U^+\mid \tau(x,U_{\BA}^{L-})\in {\BA}\},\\
\label{eq:UA-}
&U^-_{\BA}=\{y\in U^-\mid \tau(U_{\BA}^{L+},y)\in {\BA}\},\\
\label{eq:UA0}
&U^0_{\BA}=\sum_{\lambda\in P}{\BA} K_\lambda.
\end{align}
Moreover, the multiplication of $U$ induces the isomorphism
\begin{align*}
&U_{\BA}\simeq
U_{\BA}^{+}\otimes
U_{\BA}^{0}\otimes
U_{\BA}^{-},\\
&U_{\BA}^{\geqq0}\simeq
U_{\BA}^{0}\otimes
U_{\BA}^{+}
\simeq
U_{\BA}^{+}\otimes
U_{\BA}^{0},\\
&U_{\BA}^{\leqq0}\simeq
U_{\BA}^{0}\otimes
U_{\BA}^{-}
\simeq
U_{\BA}^{-}\otimes
U_{\BA}^{0}
\end{align*}
of ${\BA}$-modules.

For $i\in I$ we set
\begin{equation}
A_{i}=(q_{i}-q_{i}^{-1})E_{i},\qquad
B_{i}=(q_{i}-q_{i}^{-1})F_{i}.
\end{equation}
For $1\leqq k\leqq N$ we also set
\begin{equation}
A_{\beta_k}=(q_{\beta_k}-q_{\beta_k}^{-1})E_{\beta_k},\qquad
B_{\beta_k}=(q_{\beta_k}-q_{\beta_k}^{-1})F_{\beta_k}.
\end{equation}

Then we have the following results (see Gavarini \cite{Gav}).
\begin{lemma}
\label{lem:DP-L}
$\{A_{\beta_{N}}^{m_N}\cdots A_{\beta_{1}}^{m_1}\mid
m_1,\dots, m_N\geqq0\}$ (resp. $\{B_{\beta_{N}}^{m_N}\cdots B_{\beta_{1}}^{m_1}\mid
m_1,\dots, m_N\geqq0\}$)
is an ${\BA}$-basis of $U_{\BA}^+$ (resp. $U_{\BA}^-$).
In particular, we have $U_{\BA}^{\pm}\subset U_{\BA}^{L, \pm}$.
\end{lemma}
\begin{proposition}
\begin{itemize}
\item[\rm(i)]
$U_{\BA}^{0}$,  $U_{\BA}^+$, $U_{\BA}^-$, $U_{\BA}^{\geqq0}$,  $U_{\BA}^{\leqq0}$, $U_{\BA}$ are ${\BA}$-subalgebras of $U$.
\item[\rm(ii)]
$U_{\BA}^0$,  
$U_{\BA}^{\geqq0}$,  $U_{\BA}^{\leqq0}$, $U_{\BA}$ 
are Hopf algebras over ${\BA}$.
\end{itemize}
\end{proposition}
\begin{remark}{\rm
It follows from Lemma  \ref{lem:DP-L} that the ${\BA}$-form $U_{\BA}$ of $U$ is the same as the one considered in De Concini-Procesi \cite{DP}.}
\end{remark}

By \eqref{eq:UA+}, \eqref{eq:UA-} and Lemma \ref{lem:DP-L} the bilinear form $\tau:U^{\geqq0}\times U^{\leqq0}\to\BF$ induces 
\begin{align}
\label{eq:tau-A1}
&\tau^{\emptyset,L}_{\BA}:U_{\BA}^{\geqq0}\times U_{\BA}^{L,\leqq0}\to {\BA},\qquad
\tau_{\BA}^{L,\emptyset}:U_{\BA}^{L, \geqq0}\times U_{\BA}^{\leqq0}\to {\BA},\\
\label{eq:tau-A2}
&\tau^{\emptyset,\emptyset}_{\BA}:U_{\BA}^{\geqq0}\times U_{\BA}^{\leqq0}\to {\BA}
\end{align}
and $\sigma:U\times V\to\BF$ induces a bilinear form
\begin{equation}
\label{eq:sigma-A}
\sigma_{\BA}:U_{\BA}\times V_{\BA}\to {\BA}.
\end{equation}

\section{Specialization}
Fix $z\in\BC^\times$ and set
\[
\BA_z=\{{f}/{g}\mid
f,g\in\BQ[q,q^{-1}], g(z)\ne0\}\subset\BF.
\]
We set
\begin{equation}
U^L_z=\BC\otimes_{\BA_z} U^L_{\BA_z}, \qquad
V_z=\BC\otimes_{\BA_z} V_{\BA_z},\qquad
U_z=\BC\otimes_{\BA_z} U_{\BA_z},
\end{equation}
where the specialization ${\BA_z}\to \BC$ is given by $q\mapsto z$.
We also define 
$U_z^{L,\pm}$, $U_z^{L,0}$, $U_z^{L,\geqq0}$, $U_z^{L,\leqq0}$, $U_{z,\pm\beta}^{L,\pm}\,\,(\beta\in Q^+)$, 
$V_z^{\pm}$, $V_z^{0}$, $V_z^{\geqq0}$, $V_z^{\leqq0}$, $V_{z,\pm\beta}^{\pm}\,\,(\beta\in Q^+)$,
$U_z^{\pm}$, $U_z^{0}$, $U_z^{\geqq0}$, $U_z^{\leqq0}$, $U_{z,\pm\beta}^{\pm}\,\,(\beta\in Q^+)$
similarly.
We denote by
\begin{equation}
p^L_z:U_{\BA_z}^L\to U_z^L,\qquad
p_z:V_{\BA_z}\to V_z,\qquad
\pi_z:U_{\BA_z}\to U_z\qquad
\end{equation}
the natural homomorphisms.
The bilinear forms \eqref{eq:tau-A1}, \eqref{eq:tau-A2} for $R=\BA_z$ induce bilinear forms
\begin{align}
\label{eq:tau-zeta1}
&\tau^{\emptyset,L}_z:U_z^{\geqq0}\times U_z^{L,\leqq0}\to\BC,\qquad
\tau_z^{L,\emptyset}:U_z^{L, \geqq0}\times U_z^{\leqq0}\to\BC,\\
\label{eq:tau-zeta2}
&\tau^{\emptyset,\emptyset}_z:U_z^{\geqq0}\times U_z^{\leqq0}\to\BC,
\end{align}
and \eqref{eq:sigma-A} for $R=\BA_z$ induces a bilinear form 
\begin{equation}
\label{eq:sigma-zeta}
\sigma_z:U_z\times V_z\to\BC.
\end{equation}
Set 
\[
J_z=\{v\in V_z\mid\sigma_z(U_z,v)=\{0\}\}.
\]
\begin{lemma}
\label{lem:J}
$J_z$ is a Hopf ideal of $V_z$, and we have
$J_z=V_z^-V_z^+(J_z\cap V^0_z)$.
\end{lemma}
\begin{proof}
It easily follows from Proposition \ref{prop:invariance} that $J_z$ is a two-sided ideal.

Set 
$J_z^0=J_z\cap V^0_z$ and
$
V'_z=V_z/V_z^-V_z^+J_z^0$.
Since the multiplication of $V_z$ induces an isomorphism
$V_z\simeq V_z^-\otimes V_z^+\otimes V^0$,
we have
\[
V'_z\simeq
(V_z^-\otimes V_z^+\otimes V^0)/
(V_z^-\otimes V_z^+\otimes J_z^0)
\simeq
V_z^-\otimes V_z^+\otimes (V_z^0/J_z^0).
\]
Let
$\sigma_z':U_z\times V'_z\to\BC$ be the bilinear form induced by $\sigma_z$.
Then we see easily that $\{v\in V'_z\mid\sigma'_z(U_z,v)=\{0\}\}=\{0\}$.
Hence $J_z=V_z^-V_z^+J_z^0$.
It remains to show $\Delta(J_z)\subset J_z\otimes V_z+
V_z\otimes J_z$.
By the above argument we are reduced to showing
$\Delta(J_z^0)\subset J^0_z\otimes V^0_z+
V^0_z\otimes J^0_z$.
We can easily check this from the definition of $\sigma$.
\end{proof}
We define a Hopf algebra $\overline{V}_z$ by
\begin{equation}
\overline{V}_z=V_z/J_z.
\end{equation}
We denote by 
\begin{equation}
\overline{p}_z:V_{\BA_z}\to \overline{V}_z
\end{equation}
the canonical homomorphism.
Let
\begin{equation}
\label{eq:bar-sigma-zeta}
\overline{\sigma}_z:U_z\times \overline{V}_z\to\BC
\end{equation}
be the bilinear form induced by \eqref{eq:sigma-zeta}.
Denote the images of $V_z^0$, $V_z^\pm$, $V_z^{\geqq0}$, $V_z^{\leqq0}$ under $V_z\to\overline{V}_z$ by $\overline{V}_z^0$, $\overline{V}_z^\pm$, $\overline{V}_z^{\geqq0}$, $\overline{V}_z^{\leqq0}$ respectively.
Then the multiplication of $\overline{V}_z$ induces isomorphisms
\begin{align*}
&
\overline{V}_z
\simeq
\overline{V}_z^-\otimes \overline{V}_z^+\otimes \overline{V}_z^0,\\
&\overline{V}_z^{\geqq0}
\simeq
\overline{V}_z^+\otimes \overline{V}_z^0,\\
&\overline{V}_z^{\leqq0}
\simeq
\overline{V}_z^-\otimes \overline{V}_z^0.
\end{align*}
\begin{lemma}
\label{lem:perfect}
The bilinear form \eqref{eq:bar-sigma-zeta} is perfect in the sense that
\begin{align}
\label{eq:perfect1}
&u\in U_z, \quad\overline{\sigma}_z(u,\overline{V}_z)=\{0\}
\quad\Longrightarrow \quad u=0,\\
\label{eq:perfect2}
&v\in \overline{V}_z, \quad\overline{\sigma}_z(U_z,v)=\{0\}
\quad\Longrightarrow \quad v=0.
\end{align}
\end{lemma}
\begin{proof}
\eqref{eq:perfect2} is clear from the definition.
The proof of \eqref{eq:perfect1} is reduced to showing
\[
u\in U^0_z, \quad{\sigma}_z(u,{V}^0_z)=\{0\}
\quad\Longrightarrow \quad u=0.
\]
This can be shown by a direct computation.
Details are omitted.
\end{proof}

Set 
\begin{align*}
&I_z^0=\eta^{\geqq0}(J_z\cap V_z^0)\subset U^{L,0}_z,\\
&I_z^{\geqq0}=U_z^{L,+}I_z^0
\subset U^{L,\geqq0}_z,\qquad
I_z^{\leqq0}=U_z^{L,-}I_z^0
\subset U^{L,\leqq0}_z,\\
&I_z=U_z^{L,-}U_z^{L,+}I_z^0
\subset U^{L}_z.
\end{align*}
\begin{lemma}
$I_z^0$, 
$I_z^{\geqq0}$, 
$I_z^{\leqq0}$, 
$I_z$
are Hopf ideals of 
$U^{L,0}_z$, 
$U^{L,\geqq0}_z$, 
$U^{L,\leqq0}_z$,
$U^{L}_z$  respectively.
\end{lemma}
\begin{proof}
By Lemma \ref{lem:J} we see easily that 
$J_z^0$,
$V_z^+J_z^0$, 
$V_z^-J_z^0$ 
are Hopf ideals of $V_z^0$, $V_z^{\geqq0}$,
$V_z^{\leqq0}$ respectively.
Since $\eta^{\geqq0}|_{V_z^0}:V_z^0\to U^{L,0}_z$, 
$\eta^{\geqq0}:V_z^{\geqq0}\to U^{L,\geqq0}_z$, 
$\eta^{\leqq0}:V_z^{\leqq0}\to U^{L,\leqq0}_z$
are isomorphisms of Hopf algebras,
$I_z^0$, 
$I_z^{\geqq0}$, 
$I_z^{\leqq0}$
are Hopf ideals of 
$U^{L,0}_z$, 
$U^{L,\geqq0}_z$, 
$U^{L,\leqq0}_z$ respectively.
Then the assertion for $I_z$ follows from those for $U^{L,\geqq0}_z$ and
$U^{L,\leqq0}_z$.
\end{proof} 
We define a Hopf algebra $\overline{U}^L_z$ by
\begin{equation}
\overline{U}^L_z=U^L_z/I_z.
\end{equation}
We denote by 
\begin{equation}
\overline{p}^L_z:U^L_{\BA_z}\to \overline{U}^L_z
\end{equation}
the canonical homomorphism.
Denote the images of $U^{L,0}_z$, $U_z^{L,\pm}$, $U_z^{L,\geqq0}$, $U_z^{L,\leqq0}$ under $U^L_z\to\overline{U}^L_z$ by $\overline{U}^{L,0}_z$, $\overline{U}_z^{L,\pm}$, $\overline{U}_z^{L,\geqq0}$, $\overline{U}_z^{L,\leqq0}$ respectively.
We also denote by
\begin{equation}
\overline{\eta}_z^{\geqq0}:\overline{V}_z^{\geqq0}\to\overline{U}_z^{L,\geqq0},
\qquad
\overline{\eta}_z^{\leqq0}:\overline{V}_z^{\leqq0}\to\overline{U}_z^{L,\leqq0}
\end{equation}
the Hopf algebra isomorphisms induced by $\eta^{\geqq0}$ and $\eta^{\leqq0}$.
The bilinear forms \eqref{eq:tau-zeta1} induce
\begin{equation}
\label{eq:bar-tau-zeta1}
\overline{\tau}^{\emptyset,L}_z:U_z^{\geqq0}\times \overline{U}_z^{L,\leqq0}\to\BC,\qquad
\overline{\tau}_z^{L,\emptyset}:\overline{U}_z^{L, \geqq0}\times U_z^{\leqq0}\to\BC.
\end{equation}

\section{Specialization to $1$}

For an algebraic groups $S$ over $\BC$ with Lie algebra $\Gs$ we will identify the coordinate algebra $\BC[S]$ of $S$ with a subspace of the dual space $U(\Gs)^*$ of the enveloping algebra $U(\Gs)$ by the canonical Hopf paring 
\[
\BC[S]\otimes U(\Gs)\to \BC.
\]

We see easily that $J_1$ is generated by the elements
$p_1(Z_\lambda)\in V_1$ for $\lambda\in Q$.
From this we see easily the following.
\begin{lemma}
\label{lem:barv1}
\begin{itemize}
\item[{\rm(i)}]
We have an isomorphism $\overline{V}_1\cong U(\Gm)$ of Hopf algebras satisfying
\begin{align*}
&\overline{p}_1({X}_i)\leftrightarrow x_i,\qquad
\overline{p}_1({Y}_i)\leftrightarrow y_i,\\
&\overline{p}_1\left({\begin{bmatrix}Z_i\\m\end{bmatrix}}\right)\leftrightarrow
\begin{pmatrix}t_i\\m\end{pmatrix}:=t_i(t_i-1)\cdots(t_i-m+1)/m!.
\end{align*}
\item[\rm(ii)]
We have an isomorphism $\overline{U}^L_1\cong U(\Gg)$ of Hopf algebras satisfying
\begin{align*}
&\overline{p}_1^L({E}_i)\leftrightarrow e_i,\qquad
\overline{p}_1^L({F}_i)\leftrightarrow f_i,\\
&\overline{p}_1^L\left({\begin{bmatrix}K_i\\m\end{bmatrix}}\right)\leftrightarrow
\begin{pmatrix}h_i\\m\end{pmatrix}:=h_i(h_i-1)\cdots(h_i-m+1)/m!.
\end{align*}
\end{itemize}
\end{lemma}
In the rest of this paper we will occasionally identify $\overline{V}_1$ and $\overline{U}^L_1$ with $U(\Gm)$ and $U(\Gg)$ respectively.

In \cite{DP} De Concini-Procesi proved an isomorhphism
\begin{equation}
\label{eq:DP}
U_1\cong \BC[M]
\end{equation}
of Poisson Hopf algebras.
They established \eqref{eq:DP} by giving a correspondence between generators of both sides and proving the compatibility after a lengthy calculation.
Later Gavarini \cite{Gav} gave a more natural approach to the isomorphism \eqref{eq:DP}  using the Drinfeld paring.
Namely we have the following.
\begin{theorem}[Gavarini \cite{Gav}]
\label{thm:isom}
The bilinear form 
\[
\overline{\sigma}_1:U_1\times\overline{V}_1\to\BC
\]
induces an
 an isomorphism
\begin{equation}
\label{eq:ISOM}
\Upsilon:U_1\to\BC[M]\,(\,\subset U(\Gm)^*\simeq\overline{V}_1^*)
\end{equation} 
of Hopf algebras.
\end{theorem}
Gavarini \cite{Gav} also proved that \eqref{eq:ISOM} is an isomorphism of Poisson algebras.
This point will be discussed later in Section \ref{sec:Poisson} below.

For convenience of readers we give a more concrete description of $\Upsilon$.
Let us first give a description of $\BC[M]$ as a subspace of $U(\Gm)^*$.
The enveloping algebra $U(\Gm^\pm)$ has the direct sum decomposition
\[
U(\Gm^\pm)=\bigoplus_{\beta\in Q^+}U(\Gm^\pm)_{\pm\beta},
\]
where
\[
U(\Gm^\pm)_{\pm\beta},
=\{x\in U(\Gm^\pm)\mid [\iota(h),x]=\beta(h)x\,\,(h\in\Gh)\}
\]
for $\beta\in Q^+$.
Then we have 
\[
\BC[M^\pm]=\bigoplus_{\beta\in Q^+}(U(\Gm^\pm)_{\pm\beta})^*
\subset U(\Gm^\pm)^*.
\]
Moreover, 
we have 
\[
\BC[M^0]=\bigoplus_{\lambda\in Q}\BC\chi_\lambda\subset U(\Gh)^*
\]
where $\chi_\lambda:U(\Gm^0)\to\BC$ is the algebra homomorphism given by $\chi_\lambda(\iota(h))=\lambda(h)\,\,(h\in\Gh)$.
The isomorphism
\[
M^-\times M^+\times M^0\simeq M\qquad
((g_-,g_+,g_0)\leftrightarrow g_-g_+g_0)
\]
of algebraic varieties induced by the product of the group $M$ gives an identification
\begin{equation}
\label{eq:M-docomp}
\BC[M^+]\otimes\BC[M^-]\otimes\BC[M^0]
\simeq\BC[M]
\end{equation}
of vector spaces.
The multiplication of the algebra $U(\Gm)$ induces an identification 
\[
U(\Gm^+)\otimes U(\Gm^-)\otimes U(\Gm^0)\simeq U(\Gm).
\]
Then the canonical embedding $\BC[M]\subset U(\Gm)^*$ is given by
\begin{align*}
\BC[M]&\simeq
\BC[M^+]\otimes\BC[M^-]\otimes\BC[M^0]
\subset
U(\Gm^+)^*\otimes U(\Gm^-)^*\otimes U(\Gm^0)^*\\
&\subset
(U(\Gm^+)\otimes U(\Gm^-)\otimes U(\Gm^0))^*
=U(\Gm)^*.
\end{align*}

Let $\tilde{\Upsilon}:U_1\to U(\Gm)^*(\simeq\overline{V}_1^*)$ be the linear map induced by $\overline{\sigma}_1$.
Then we have
\begin{multline*}
\tilde{\Upsilon}(u_+u_0S(u_-))(v_-v_+v_0)
=\Upsilon^+(u_+)(v_-)\cdot \Upsilon^0(u_0)(v_0)\cdot
\Upsilon^-(u_-)(v_+)\\
(u_\pm\in U_1^\pm, u_0\in U_1^0, 
v_\pm\in U(\Gm^\pm), v_0\in U(\Gm^0)).
\end{multline*}
where $\Upsilon^\pm:U_1^\pm\to U(\Gm^\mp)^*$ and $\Upsilon^0:U_1^0\to U(\Gm^0)^*$ are given by
\begin{align*}
\Upsilon^+(u_+)(v_-)&=\overline{\tau}_1^{\emptyset,L}(u_+,\overline{\eta}_1^{\leqq0}(v_-)),\\
\Upsilon^-(u_-)(v_+)&=\overline{\tau}_1^{L,\emptyset}(\overline{\eta}_1^{\geqq0}(v_+),u_-),\\
\Upsilon^0(\pi_1({K}_\lambda))&=\chi_{\lambda}.
\end{align*}

For $i\in I$ we define $a_i\in \BC[M^-]\subset U(\Gm^-)^*$,  $b_i\in \BC[M^+]\subset U(\Gm^+)^*$ by
\begin{align*}
&\langle a_i, U(\Gm^-)_{-\beta}\rangle=0\quad
(\beta\ne\alpha_i),\qquad
\langle a_i,y_i\rangle=-1,\\
&\langle b_i, U(\Gm^+)_{\beta}\rangle=0\quad
(\beta\ne\alpha_i),\qquad
\langle b_i,x_i\rangle=1.
\end{align*}
We identify $\BC[M^\pm], \BC[M^0]$ with subalgebras of $\BC[M]$ via \eqref{eq:M-docomp}, and regard $a_i, b_i, \chi_\lambda \,\,(i\in I, \lambda\in Q)$ as elements of $\BC[M]$.
We see easily the following.
\begin{lemma}
Under the identification \eqref{eq:ISOM} of Theorem \ref{thm:isom} we have
\[
\pi_1(A_{i})\leftrightarrow a_i,\quad
\pi_1(B_{i})\leftrightarrow b_i\chi_{-\alpha_i},\quad
\pi_1(K_{\lambda})\leftrightarrow \chi_{\lambda}\qquad
(i\in I, \lambda\in Q).
\]
\end{lemma}

\section{Specialization to roots of $1$}
\label{section:roots}
We fix a positive odd integer $\ell$.
We assume that $\ell$ is prime to 3 if $\Gg$ is of type $G_2$.
We denote by $\zeta\in\BC^\times$ a fixed primitive $\ell$-th root of 1.
\begin{remark}
{\rm
Denote by $U^{DK}_{\BQ[q,q^{-1}]}$ the De Concini-Kac $\BQ[q,q^{-1}]$-form of $U$
(see \cite{DK}). 
Namely $U^{DK}_{\BQ[q,q^{-1}]}$ is the ${\BQ[q,q^{-1}]}$-subalgebra of $U$ generated by
$\{K_i^{\pm1}, E_i, F_i\mid i\in I\}$.
Then we have $U_\zeta\simeq\BC\otimes_{\BQ[q,q^{-1}]} U^{DK}_{\BQ[q,q^{-1}]}$, where the specialization ${\BQ[q,q^{-1}]}\to\BC$ is given by $q\mapsto\zeta$.
}
\end{remark}
We denote by 
$\tilde{\xi}^L:U^L_\zeta\to U^L_1$
Lusztig's Frobenius morphism (see \cite{L2}).
Namely, $\tilde{\xi}^L$ is an algebra homomorphism given by
\begin{align}
\tilde{\xi}^L(\label{eq:TXL1}
{p}_\zeta^L({E_i^{(n)}}))
&=
\begin{cases}
p_1^L({E_i^{(n/\ell)}})\quad&(\ell\,|\,n)\\
0&(\ell\not|\,n),
\end{cases}\\
\label{eq:TXL2}
\tilde{\xi}^L{p}_\zeta^L({F_i^{(n)}}))
&=
\begin{cases}
p_1^L({F_i^{(n/\ell)}})\quad&(\ell\,|\,n)\\
0&(\ell\not|\,n),
\end{cases}
\\
\label{eq:TXL3}
\tilde{\xi}^L{p}_\zeta^L
\left(
{\begin{bmatrix}{K_i}\\{m}\end{bmatrix}}
\right))
&=
\begin{cases}
p_1^L
\left(
{\begin{bmatrix}{K_i}\\{m/\ell}\end{bmatrix}}
\right)
\quad&(\ell\,|\,m)\\
\,\,\,0&(\ell\not|\,m),
\end{cases}
\\
\label{eq:TXL4}
\tilde{\xi}^L{p}_\zeta^L({K}_\lambda))
&=
p_1^L({K}_\lambda)\quad(\lambda\in Q).
\end{align}
It is a Hopf algebra homomorphism.
\begin{lemma}
\label{lem:TXL}
We have $\tilde{\xi}^L(I_\zeta)\subset I_1$.
\end{lemma}
\begin{proof}
It is sufficient to show
$\tilde{\xi}(I_\zeta^0)\subset I_1^0$.
For $z\in\BC^\times$, $m=(m_i)_{i\in I}\in\BZ_{\geqq0}^I$, and ${\varepsilon}=(\varepsilon_i)_{i\in I}\in\{0, 1\}^I$
set
\[
K_{m,\varepsilon}(z)=
p_z^L
\left(
\prod_{i\in I}\left(
K_i^{\varepsilon_i}
\begin{bmatrix}{K_i}\\{m_i}\end{bmatrix}
\right)
\right)
\in U_z^{L,0}.
\]
Then any element $u$ of $U_z^{L,0}$ is uniquely written as a finite sum
\[
u=\sum_{m,\varepsilon}c_{m,\varepsilon}K_{m,\varepsilon}(z)
\qquad(c_{m,\varepsilon}\in\BC).
\]
Then we have $u\in I^0$ if an d only if
\[
\left.\sum_{m,\varepsilon}c_{m,\varepsilon}q_i^{(\lambda,\alpha_i^\vee)}
\begin{bmatrix}{(\lambda,\alpha_i^\vee)}\\{m_i}\end{bmatrix}_{q_i}
\right|_{q=z}=0
\qquad(\forall\lambda\in Q).
\]
Hence it is sufficient to show that
\begin{equation}
\label{eq:J1}
\left.\sum_{m,\varepsilon}c_{m,\varepsilon}q_i^{(\lambda,\alpha_i^\vee)}
\begin{bmatrix}{(\lambda,\alpha_i^\vee)}\\{m_i}\end{bmatrix}_{q_i}
\right|_{q=\zeta}=0
\qquad(\forall\lambda\in Q)
\end{equation}
implies
\begin{equation}
\label{eq:J2}
\sum_{m,\varepsilon}c_{\ell m,\varepsilon}
\begin{pmatrix}{(\mu,\alpha_i^\vee)}\\{\ell m_i}\end{pmatrix}=0
\qquad(\forall\mu\in Q).
\end{equation}
Indeed \eqref{eq:J2} follows by setting $\lambda=\ell\mu$ in \eqref{eq:J1}.
\end{proof}
We denote by
\begin{equation}
{\xi}^L:\overline{U}^L_\zeta\to \overline{U}^L_1
\end{equation}
the Hopf algebra homomorphism
induced by $\tilde{\xi}^L$.
\begin{lemma}
There exists a Hopf algebra homomorphism
\begin{equation}
\xi:\overline{V}_\zeta\to\overline{V}_1
\end{equation}
satisfying 
\begin{align*}
\xi(\overline{p}_\zeta({X_i^{(n)}}))
&=
\begin{cases}
\overline{p}_1({X_i^{(n/\ell)}})\quad&(\ell\,|\,n)\\
0&(\ell\not|\,n),
\end{cases}
\\
\xi(\overline{p}_\zeta({Y_i^{(n)}}))
&=
\begin{cases}
\overline{p}_1({Y_i^{(n/\ell)}})\quad&(\ell\,|\,n)\\
0&(\ell\not|\,n),
\end{cases}
\\
\xi\left(
\overline{p}_\zeta
\left(
{\begin{bmatrix}{Z_i}\\{m}\end{bmatrix}}
\right)
\right)
&=
\begin{cases}
\overline{p}_1
\left({\begin{bmatrix}{Z_i}\\{m/\ell}\end{bmatrix}}
\right)
\quad&(\ell\,|\,m)\\
\,\,\,0&(\ell\not|\,m),
\end{cases}
\\
\xi((\overline{p}_\zeta({Z}_\lambda))
&=
1\quad(\lambda\in Q).
\end{align*}
\end{lemma}
\begin{proof}
By the isomorphisms $\overline{V}^{\geqq0}_z\simeq \overline{U}^{L,\geqq0}_z$, $\overline{V}^{\leqq0}_z\simeq \overline{U}^{L,\leqq0}_z$ for $z\in\BC^\times$ induced by $\eta^{\geqq0}$, $\eta^{\leqq0}$ we obtain Hopf algebra homomorphisms
$\xi^{\geqq0}:\overline{V}_\zeta^{\geqq0}\to \overline{V}_1^{\geqq0}$ and
$\xi^{\leqq0}:\overline{V}_\zeta^{\leqq0}\to \overline{V}_1^{\leqq0}$
corresponding to
$\xi^{L}|_{\overline{U}_\zeta^{L,\geqq0}}$ and $\xi^{L}|_{\overline{U}_\zeta^{L,\leqq0}}$ respectively.
By $[\overline{V}^+_z, \overline{V}^-_z]=0$ we obtain the desired Hopf algebra homomorphism $\xi:\overline{V}_\zeta\to \overline{V}_1$ satisfying
$\xi|_{\overline{V}_\zeta^{\geqq0}}=\xi^{\geqq0}$ and 
$\xi|_{\overline{V}_\zeta^{\leqq0}}=\xi^{\leqq0}$.
\end{proof}
By \cite{L2} and the construction of $\xi$ we have the following.
\begin{lemma}
We have
\begin{align*}
\xi(\overline{p}_\zeta({X_{\beta_k}^{(n)}}))&=
\begin{cases}
\overline{p}_1({X_{\beta_k}^{(n/\ell)}})\quad&(\ell\,|\,n)\\
0&(\ell\not|\,n),
\end{cases}
\\
\xi(\overline{p}_\zeta({Y_{\beta_k}^{(n)}}))&=
\begin{cases}
\overline{p}_1({Y_{\beta_k}^{(n/\ell)}})\quad&(\ell\,|\,n)\\
0&(\ell\not|\,n).
\end{cases}
\end{align*}
\end{lemma}
\begin{proposition}
\label{prop:UU}
There exists a unique linear map 
\begin{equation}
{}^t\xi:U_1\to U_\zeta
\end{equation}
satisfying
\begin{equation}
\label{eq:tx}
\overline{\sigma}_\zeta({}^t\xi(u),v)=
\overline{\sigma}_1(u,\xi(v))\qquad
(u\in U_1,\,\, v\in \overline{V}_\zeta).
\end{equation}
It is an injective Hopf algebra homomorphism whose image is contained in the center of $U_\zeta$.
\end{proposition}
\begin{proof}
By a direct computation 
the linear map ${}^t\xi:U_1\to U_\zeta$ defined by
\begin{multline}
\label{eq:tx2}
{}^t\xi(\pi_1(A_{\beta_N}^{r_N}\cdots A_{\beta_1}^{r_1}
K_\lambda S(B_{\beta_N}^{s_N}\cdots B_{\beta_1}^{s_1})))\\
=
\pi_\zeta(A_{\beta_N}^{\ell r_N}\cdots A_{\beta_1}^{\ell r_1}
K_{\ell \lambda} S(B_{\beta_N}^{\ell s_N}\cdots B_{\beta_1}^{\ell s_1}))
\end{multline}
satisfies \eqref{eq:tx}.
The uniqueness and the injectivity of ${}^t\xi$ follow from Lemma \ref{lem:perfect}.
${}^t\xi$ is a homomorphism of coalgebras by Proposition \ref{prop:invariance}.

Let us show that ${}^t\xi(u)$ is a central element for any $u\in U_1$.
We may assume that $u\in U_1^{0}$ or $u\in U_1^+$ or $u\in S(U_1^-)$.
If $u\in U_1^{0}$, then  ${}^t\xi(u)$ is a central element since it is a linear combination of the elements of the form $K_{\ell\lambda}\,\,(\lambda\in Q)$.
Assume $u\in U_1^+$. 
Let us show 
\begin{equation}
\label{eq:UU1}
{}^t\xi(u)x=x{}^t\xi(u)\qquad(x\in U_\zeta^+).
\end{equation}
It is sufficient to show 
$\overline{\sigma}_\zeta({}^t\xi(u)x,y)
=\overline{\sigma}_\zeta(x{}^t\xi(u),y)$ 
for any $y\in \overline{V}_\zeta^-$.
By
\begin{align*}
&\overline{\sigma}_\zeta({}^t\xi(u)x,y)
=(\overline{\sigma}_\zeta\otimes \overline{\sigma}_\zeta)
(x\otimes {}^t\xi(u),\Delta(y))
=
(\overline{\sigma}_\zeta
\otimes
\overline{\sigma}_1)
(x\otimes u,(1\otimes\xi)(\Delta(y))),\\
&\overline{\sigma}_\zeta(x{}^t\xi(u),y)
=(\overline{\sigma}_\zeta\otimes \overline{\sigma}_\zeta)
(x\otimes {}^t\xi(u),\Delta'(y))
=
(\overline{\sigma}_\zeta
\otimes
\overline{\sigma}_1)
(x\otimes u,(1\otimes\xi)(\Delta'(y))),
\end{align*}
it is sufficient to show $(1\otimes\xi)(\Delta(y))=(1\otimes\xi)(\Delta'(y))$ for any $y\in \overline{V}_\zeta^-$.
Here $\Delta'$ is the opposite comultiplication.
We may assume $y=\overline{p}_\zeta({Y_i^{(n)}})$.
Then we have
\[
\Delta(\overline{p}_\zeta({Y_i^{(n)}}))
=\sum_{r=0}^n
\zeta_i^{r(n-r)}\overline{p}_\zeta({Y_i^{(r)}})\otimes \overline{p}_\zeta({Y_i^{(n-r)}Z_i^{r}}),
\]
where $\zeta_i=\zeta^{(\alpha_i,\alpha_i)/2}$, and hence
\[
(1\otimes\xi)(\Delta(\overline{p}_\zeta({Y_i^{(n)}}))))
=\sum_{r\geqq0, \ell r\leqq n}
\zeta_i^{nr}
\overline{p}_\zeta({Y_i^{(n-\ell r)}})\otimes \overline{p}_\zeta({Y_i^{(r)}})
=(1\otimes\xi)(\Delta'(\overline{Y_i^{(n)}})))
\]
(note that $\overline{p}_\zeta({Z_i^\ell})=1$).
\eqref{eq:UU1} is proved.
By Proposition \ref{prop:DPPBW} we have
\begin{equation}
\label{eq:tau2}
\overline{\tau}^{\emptyset,\emptyset}_\zeta({}^t\xi(x),y)=\varepsilon(x)\varepsilon(y)\qquad(x\in U_1^{\geqq0}, y\in U_\zeta^{\leqq0}).
\end{equation}
Hence by Lemma \ref{lem:Drinfeld paring} we obtain
\begin{align*}
{}^t\xi(u)y=&
\sum_{(u)_2, (y)_2}
\overline{\tau}^{\emptyset,\emptyset}_\zeta({}^t\xi(u_{(0)}),y_{(0)})
\overline{\tau}^{\emptyset,\emptyset}_\zeta({}^t\xi(u_{(2)}),S(y_{(2)}))
y_{(1)}{}^t\xi(u_{(1)})\\
=&
\sum_{(u)_2, (y)_2}
\varepsilon(u_{(0)})\varepsilon(y_{(0)})
\varepsilon(u_{(2)})\varepsilon(y_{(2)})
y_{(1)}{}^t\xi(u_{(1)})=y{}^t\xi(u)
\end{align*}
for any $y\in U_\zeta^{\leqq0}$.
Therefore, ${}^t\xi(u)$ is a central element for any $u\in U_1^+$.
Similarly, we can show that ${}^t\xi(u)$ is a central element for any $u\in S(U_1^-)$.
We have shown that the image of ${}^t\xi$ is contained in the center of $U_\zeta$.
It also follows from this and \eqref{eq:tx2} that ${}^t\xi$ is an algebra homomorphism.
\end{proof}
\begin{remark}
{\rm
Some of the arguments in our proof of Proposition \ref{prop:UU} is similar to those for the dual statement in Gavarini \cite[Theorem 7.9]{Gav}.
}
\end{remark}
We set
\begin{equation}
Z_\zeta=\Image({}^t\xi).
\end{equation}
By Proposition \ref{prop:UU} it is a Hopf subalgebra of $U_\zeta$ contained in the center.

\section{Poisson structures}
\label{sec:Poisson}
By Theorem \ref{thm:isom} and Proposition \ref{prop:UU} we have isomorphisms
\begin{equation}
\label{ISOMORPHISM}
Z_\zeta\simeq U_1\simeq\BC[M]
\end{equation}
of Hopf algebras.
They are in fact isomorphisms of Poisson Hopf algebras with respect to certain canonical Poisson structures (De Concini-Procesi \cite{DP}).
In this section we will give an account of those Poisson structures.

We first recall standard facts on Poisson structures (see e.g.\ \cite{DP}).
A commutative associative algebra $\CR$ over a field $\BK$ equipped with a bilinear map
\[
\{\,,\,\}:\CR\times \CR\to \CR
\]
is called a Poisson algebra if it satisfies
\begin{itemize}
\item[(a)]
$\{a,a\}=0\quad(a\in \CR)$,
\item[(b)]
$\{a,\{b,c\}\}+\{b,\{c,a\}\}+\{c,\{a,b\}\}=0\quad(a, b, c\in \CR)$,
\item[(c)]
$\{a,bc\}=b\{a,c\}+\{a,b\}c\quad(a, b, c\in \CR)$.
\end{itemize}
A map $F:\CR\to \CR'$ between Poisson algebras $\CR$, $\CR'$ is called a homomorphism of Poisson algebras if it is a homomorphism of associative algebras and satisfies
$F(\{a_1,a_2\})=\{F(a_1),F(a_2)\}$ for any $a_1, a_2\in \CR$.
The tensor product $\CR\otimes_\BK \CR'$ of two Poisson algebras $\CR$, $\CR'$ over $\BK$ is equipped with a canonical Poisson algebra structure given by 
\begin{align*}
&(a_1\otimes b_1)(a_2\otimes b_2)=a_1a_2\otimes b_1b_2,\\
&\{a_1\otimes b_1,a_2\otimes b_2\}=\{a_1,a_2\}\otimes b_1b_2+a_1a_2\otimes \{b_1,b_2\}
\end{align*}
for $a_1, a_2\in \CR$, $b_2, b_2\in \CR'$.
A commutative Hopf algebra $\CR$  over a field $\BK$ equipped with a bilinear map
\[
\{\,,\,\}:\CR\times \CR\to \CR
\]
is called a Poisson Hopf algebra if it is a Poisson algebra and the comultiplication $\Delta:\CR\to \CR\otimes_\BK \CR$ is a homomorphism of Poisson algebras (in this case the counit $\varepsilon:\CR\to\BK$ and the antipode $S:\CR\to \CR$ become automatically a homomorphism and an anti-homomorphism of Poisson algebras respectively).

A smooth affine algebraic variety $X$ over $\BC$ is called a Poisson variety if we are given a bilinear map
\[
\{\,,\,\}:\BC[X]\times\BC[X]\to\BC[X]
\]
so that $\BC[X]$ is a Poisson algebra.
In this case $\{f, g\}(x)$ for $f, g\in\BC[X]$ and $x\in X$ depends only on $df_x, dg_x$, and hence we have $\delta\in\Gamma(X,\bigwedge^2\Theta_X)$ 
such that 
\[
\{f, g\}(x)=\delta_x(df_x, dg_x),
\]
where $\Theta_X$ denotes the sheaf of vector fields.
We call $\delta$ the Poisson tensor of the Poisson variety $X$.

A linear algebraic group $S$ over $\BC$ is called a Poisson algebraic group if we are given a bilinear map
\[
\{\,,\,\}:\BC[S]\times\BC[S]\to\BC[S]
\]
so that $\BC[S]$ is a Poisson Hopf algebra.
Let $\delta$ be the Poisson tensor of $S$ as a Poisson variety, and define $\varepsilon:S\to \bigwedge^2\Gs$ by 
$(d\ell_g)(\varepsilon(g))=\delta_g$ for $g\in S$.
Here, $\Gs$ is the Lie algebra of $S$ which is identified with the tangent space $T_eS$ at the identity element $e\in S$, and $\ell_g:S\to S$ is given by $x\mapsto gx$.
By differentiating $\varepsilon$ at $e$ we obtain a linear map $\Gs\to\bigwedge^2\Gs$.
It induces an alternating bilinear map $[\,,\,]:\Gs^*\times\Gs^*\to\Gs^*$.
Then this $[\,,\,]$ gives a Lie algebra structure on $\Gs^*$.
Moreover, the following bracket product gives a Lie algebra structure on $\Gs\oplus\Gs^*$:
\[
[(a,\varphi),(b,\psi)]
=([a,b]+\varphi b-\psi a, a\psi-b\varphi+[\varphi,\psi]).
\]
Here, $\Gs\times \Gs^*\ni(a,\varphi)\to a\varphi\in\Gs^*$ and 
$\Gs^*\times \Gs\ni(\varphi, a)\to \varphi a\in\Gs$ are the coadjoint actions of $\Gs$ and $\Gs^*$ on $\Gs^*$ and $\Gs$ respectively.
In other words $(\Gs\oplus\Gs^*,\Gs,\Gs^*)$ is a Manin triple with respect to the symmetric bilinear form on $\Gs\oplus\Gs^*$ given by 
$((a,\varphi),(b,\psi))=\varphi(b)+\psi(a)$.
We say that $(\Ga,\Gb,\Gc)$ is a Manin triple with respect to a symmetric bilinear form $(\,,\,)$ on $\Ga$ if 
\begin{itemize}
\item[(a)]
$\Ga$ is a Lie algebra,
\item[(b)]
$(\,,\,)$ is $\Ga$-invariant and non-degenerate,
\item[(c)]
$\Gb$ and $\Gc$ are subalgebras of $\Ga$ such that $\Ga=\Gb\oplus\Gc$ as a vector spaces,
\item[(d)]
$(\Gb,\Gb)=(\Gc,\Gc)=\{0\}$.
\end{itemize}
Conversely, if $(\Ga,\Gb,\Gc)$ is a Manin triple and $B$ is a linear algebraic group with Lie algebra $\Gb$, then we can associate a natural Poisson Hopf algebra structure on $\BC[B]$ by reversing the above process.

Now let us return to our original setting.
Note that $\Gm$ is a subalgebra of $\Gg\oplus\Gg$.
Set 
\begin{equation}
\Gk=\{(x,x)\mid x\in\Gg\}\subset\Gg\oplus\Gg.
\end{equation} 
We have a natural isomorphism
\begin{equation}
\theta:\Gg\to\Gk\qquad(\theta(x)=(x,x))
\end{equation}
of Lie algebras.
It is easily seen that $(\Gg\oplus\Gg,\Gm,\Gk)$ is a Manin triple with respect to the symmetric bilinear form $\tilde{\kappa}$ on $\Gg\oplus\Gg$ given by 
\begin{equation}
\tilde{\kappa}((x_1,y_1),(x_2,y_2))=\kappa(x_1,x_2)-\kappa(y_1,y_2)\qquad
(x_1, x_2, y_1, y_2\in\Gg),
\end{equation}
where 
\begin{equation}
\kappa:\Gg\times \Gg\to\BC
\end{equation}
is the $\Gg$-invariant symmetric bilinear form which induces the symmetric bilinear form \eqref{eq:SBMh*} on $\Gh^*$.
It follows that $\BC[M]$ is endowed with a natural Poisson Hopf algebra structure.
\begin{lemma}
\label{lem:P-gen}
$\BC[M]$ is generated by 
$\{a_i, b_i, \chi_\lambda\mid i\in I, \lambda\in Q\}$ as a Poisson algebra.
\end{lemma}
\begin{proof}
See De Concini-Procesi \cite[Section 14.5]{DP}.
\end{proof}

On the other hand we have a natural Poisson Hopf algebra structure on $U_1$ given by 
\begin{equation}
\label{eq:PoissonU1}
\{\pi_1({a}),\pi_1({b})\}=\pi_1
\left(
[a,b]/(q-q^{-1})
\right)\qquad
(a, b\in U_{\BA_1})
\end{equation}
(see De Concini-Procesi \cite{DP}).

The definition of the Poisson structure on $Z_\zeta$ is more subtle.
Let $C_\zeta$ be the center of $U_\zeta$.
We have a Poisson algebra structure on $C_\zeta$ given by
\begin{multline}
\label{eq:PoissonUzeta}
\{\pi_\zeta({a}),\pi_\zeta({b})\}=\pi_\zeta
\left(
[a,b]/\ell(q^\ell-q^{-\ell})\right)\\
(a, b\in U_{\BA_\zeta},\,\, \pi_\zeta({a}), \pi_\zeta({b})\in C_\zeta).\qquad
\end{multline}
If $Z_\zeta$ is closed under the Poisson bracket \eqref{eq:PoissonUzeta}, then this gives a Poisson Hopf algebra structure on $Z_\zeta$ (see De Concini-Procesi \cite{DP}).

\begin{theorem}[De Concini-Procesi \cite{DP}]
\label{thm:Poisson}
$Z_\zeta$ is closed under the Poisson bracket \eqref{eq:PoissonUzeta}.
Moreover, 
the isomorphisms in \eqref{ISOMORPHISM}
preserve Poisson structures.
\end{theorem}
Gavarini \cite{Gav} also gave a natural proof of the fact that the isomorphism $U_1\simeq \BC[M]$ in \eqref{ISOMORPHISM} preserves the Poisson structures  using his definition  of the isomorphism in terms of the Drinfeld paring.
In fact, he gave a proof of the statement dual to it concernig Poisson coalgebra structure of the dual objects, and deduced the above statement from it.
In the rest of this section we give a direct proof of this statement.

\begin{lemma}
\label{lem:PU+}
Let $i\in I$, $\gamma\in Q^+$  and $b\in U^-_{1,-\gamma}$.
Write
\[
\Delta(b)=
\sum_
{\gamma_1, \gamma_2\in Q^+,\gamma_1+\gamma_2=\gamma}
b_{\gamma_1,\gamma_2}(1\otimes \pi_1(K_{-\gamma_1}))
\qquad
(b_{\gamma_1,\gamma_2}\in U^-_{1,-\gamma_1}\otimes U^-_{1,-\gamma_2}),
\]
and define  $b', b''\in U_{1,-(\gamma-\alpha_i)}^-$ by
\begin{equation}
\label{eq:PU+}
b_{\alpha_i,\gamma-\alpha_i}=\pi_1(B_i)\otimes b',\qquad
b_{\gamma-\alpha_i,\alpha_i}=b''\otimes \pi_1(B_i).
\end{equation}
Then we have
\begin{align}
\label{lem:PU+1}
&\{\pi_1(A_i),b\}
=\frac{(\alpha_i,\alpha_i)}2(b''\pi_1(K_i)-b'\pi_1(K_i^{-1})),\\
\label{lem:PU+2}
&\{\pi_\zeta(A_i^\ell),{}^t\xi(b)\}
=\frac{(\alpha_i,\alpha_i)}2({}^t\xi(b'')\pi_\zeta(K_i^\ell)-{}^t\xi(b')\pi_\zeta(K_i^{-\ell})).
\end{align}
\end{lemma}
\begin{proof}
Note that \eqref{lem:PU+1} can be regarded as a special case of \eqref{lem:PU+2} when $\ell=1$.
Hence we will only prove \eqref{lem:PU+2}.
We can write 
\[
\Delta_2(A_i^\ell)
=A_i^\ell\otimes1\otimes1
+K_i^\ell\otimes A_i^\ell\otimes1
+K_i^\ell\otimes K_i^\ell\otimes A_i^\ell
+(q^\ell-q^{-\ell})\sum_jX_j\otimes X'_j\otimes X''_j
\]
for some
$X_j, X'_j, X''_j\in U_{\BA_\zeta}^{\geqq0}$.
By $(\varepsilon\otimes 1\otimes\varepsilon)\Delta_2(A_i^\ell)=A_i^\ell$ we have
\begin{equation}
\label{eq:eee0}
\sum_j\varepsilon(X_j)\varepsilon(X''_j)X'_j=0.
\end{equation}
Take $B\in U_{\BA,-\ell\gamma}^-$ such that $\pi_\zeta(B)={}^t\xi(b)$.
Then we have 
\begin{align*}
&A_i^\ell B\\
=&\sum_{(B)_2}
\tau(A_i^\ell,B_{(0)})\tau(1,SB_{(2)})B_{(1)}
+
\sum_{(B)_2}
\tau(K_i^\ell,B_{(0)})\tau(1,SB_{(2)})B_{(1)}A_i^\ell\\
&\qquad
+
\sum_{(B)_2}
\tau(K_i^\ell,B_{(0)})\tau(A_i^\ell,SB_{(2)})B_{(1)}K_i^\ell\\
&\qquad
+(q^{\ell}-q^{-\ell})\sum_{(B)_2,j}
\tau(X_j,B_{(0)})\tau(X_j'',SB_{(2)})B_{(1)}X'_j\\
=&
\sum_{(B)}
\tau(A_i^\ell,B_{(0)})B_{(1)}
+
BA_i^\ell
+
\sum_{(B)}
\tau(A_i^\ell,SB_{(1)})B_{(0)}K_i^\ell\\
&\qquad
+(q^{\ell}-q^{-\ell})\sum_{(B)_2,j}
\tau(X_j,B_{(0)})\tau(X_j'',SB_{(2)})B_{(1)}X'_j,
\end{align*}
and hence
\begin{align*}
&\ell\{\pi_\zeta(A_i^\ell),{}^t\xi(b)\}\\
=&
\pi_\zeta
\left(
\left.
\left(
\sum_{(B)}
\tau(A_i^\ell,B_{(0)})B_{(1)}
+
\sum_{(B)}
\tau(A_i^\ell,SB_{(1)})B_{(0)}K_i^\ell
)\right)
\right/(q^\ell-q^{-\ell})
\right)\\
\qquad
&+
\pi_\zeta
\left(
\sum_{(B)_2,j}
\tau(X_j,B_{(0)})\tau(X_j'',SB_{(2)})B_{(1)}X'_j
\right).
\end{align*}
Note
\begin{align*}
&\sum_{(B)_2,j}
\pi_\zeta(B_{(0)})\otimes
\pi_\zeta(B_{(1)})\otimes
\pi_\zeta(B_{(2)})
=(\pi_\zeta\otimes\pi_\zeta\otimes\pi_\zeta)(\Delta_2(B))\\
=&
\Delta_2(\pi_\zeta(B))
=\Delta_2({}^t\xi(b))
=({}^t\xi\otimes{}^t\xi\otimes{}^t\xi)(\Delta_2(b)).
\end{align*}
Hence by \eqref{eq:tau2} and \eqref{eq:eee0} we have
\begin{align*}
&\pi_\zeta
\left(
\sum_{(B)_2,j}
\tau(X_j,B_{(0)})\tau(X_j'',SB_{(2)})B_{(1)}X'_j
\right)\\
=&
\sum_{(B)_2,j}
\overline{\tau}^{\emptyset,\emptyset}_\zeta(\pi_\zeta(X_j),\pi_\zeta(B_{(0)}))\overline{\tau}^{\emptyset,\emptyset}_\zeta(\pi_\zeta(X_j''),S(\pi_\zeta(B_{(2)}))\pi_\zeta(B_{(1)})\pi_\zeta(X'_j)\\
=&
\sum_{(b)_2,j}
\overline{\tau}^{\emptyset,\emptyset}_\zeta(\pi_\zeta(X_j),{}^t\xi(b_{(0)}))
\overline{\tau}^{\emptyset,\emptyset}_\zeta(\pi_\zeta(X_j''),{}^t\xi(Sb_{(2)})){}^t\xi(b_{(1)})\pi_\zeta(X'_j)\\
=&
{}^t\xi(b)
\pi_\zeta
\left(
\sum_{j}
\varepsilon(X_j)\varepsilon(X''_j)X'_j
\right)=0.
\end{align*}
Write
\[
\Delta(B)=
\sum_
{\gamma_1, \gamma_2\in Q^+,\gamma_1+\gamma_2=\ell\gamma}
B_{\gamma_1,\gamma_2}(1\otimes K_{-\gamma_1})
\qquad
(B_{\gamma_1,\gamma_2}\in U^-_{\BA,-\gamma_1}\otimes U^-_{\BA,-\gamma_2}),
\]
and define  $B', B''\in U_{\BA,-\ell(\gamma-\alpha_i)}^-$ by
\[
B_{\ell\alpha_i,\ell(\gamma-\alpha_i)}=B_i^\ell\otimes B',\qquad
B_{\ell(\gamma-\alpha_i),\ell\alpha_i}=B''\otimes B_i^\ell.
\]
Then by
\[
({}^t\xi\otimes{}^t\xi)(\Delta(b))
=\Delta({}^t\xi(b))
=\Delta(\pi_\zeta(B))
=(\pi_\zeta\otimes\pi_\zeta)(\Delta(B))
\]
we have $\pi_\zeta(B')={}^t\xi(b')$ and $\pi_\zeta(B'')={}^t\xi(b'')$.
Hence we obtain
\begin{align*}
&\ell\{\pi_\zeta(A_i^\ell),{}^y\xi(b)\}\\
=&
\pi_\zeta
\left(
\left.
\left(
\tau(A_i^\ell,B_i^\ell)B'K_i^{-\ell}+
\tau(A_i^\ell,S(B_i^\ell K_{-\ell(\gamma-\alpha_i)}))B''K_i^\ell
\right)
\right/(q^\ell-q^{-\ell})
\right)\\
=&
\ell(\alpha_i,\alpha_i)
(\pi_\zeta(B'')\pi_\zeta(K_i^\ell)
-\pi_\zeta(B')\pi_\zeta(K_i^{-\ell})
)/2
\\
=&
\ell(\alpha_i,\alpha_i)
({}^t\xi(b'')\pi_\zeta(K_i^\ell)
-{}^t\xi(b')\pi_\zeta(K_i^{-\ell})
)/2.
\end{align*}
\end{proof}
For $i\in I$ we set
\begin{equation}
\zeta_i=\zeta^{(\alpha_i,\alpha_i)/2}.
\end{equation}
For $F\in U^{L,-}_{\BA_\zeta,-\gamma}\,\,(\gamma\in Q^+)$ define 
$\varphi^i_{r,s}(F)\in U^{L,-}_{\BA_\zeta,-(\gamma-(r+s)\alpha_i)}\,\,(r,s\geqq0,\,\,i\in I)$ by
\begin{multline}
\label{eq:phi-rs}
\Delta_2(F)\in
\sum_{r,s}
F_i^{(r)}\otimes \varphi^i_{r,s}(F)K_i^{-r}\otimes F_i^{(s)}K_{-\gamma+s\alpha_i}\\
+
\bigoplus
_{(\gamma_1,\gamma_2,\gamma_3)\in \Xi}
U^{L,-}_{\BA_\zeta,-\gamma_1}\otimes
U^{L,-}_{\BA_\zeta,-\gamma_2}K_{-\gamma_1}\otimes
U^{L,-}_{\BA_\zeta,-\gamma_3}K_{-\gamma_1-\gamma_2},
\end{multline}
where $\Xi$ consists of $(\gamma_1,\gamma_2,\gamma_3)\in(Q^+)^3$ such that $\gamma_1+\gamma_2+\gamma_3=\gamma$ and $(\gamma_1,\gamma_3)\notin\BZ_{\geqq0}\alpha_i\times\BZ_{\geqq0}\alpha_i$.
\begin{lemma}
\label{lem:varphi}
Let $i\in I$.
For $F\in U^{L,-}_{\BA_\zeta,-\gamma}\,\,(\gamma\in Q^+)$ we have
\[
\xi^L(\overline{p}^L_\zeta(\varphi^i_{r,s}(F)))=\zeta_i^{rs}
\xi^L(\overline{p}^L_\zeta(\varphi^i_{r+s,0}(F))).
\]
\end{lemma}
\begin{proof}
For $X\in U^+_{\BA_\zeta}$ we have
\begin{align*}
\tau(A_i^sXA_i^r,F)
=&(\tau\otimes\tau\otimes\tau)(A_i^r\otimes X\otimes A_i^s,\Delta_2(F))\\
=&\tau(A_i^r,F_i^{(r)})\tau(X,\varphi^i_{r,s}(F))\tau(A_i^s,F_i^{(s)})\\
=&(-1)^{r+s}q_i^{(r(r-1)+s(s-1))/2}\tau(X,\varphi^i_{r,s}(F)),
\end{align*}
and hence
\[
\overline{\tau}^{\emptyset,L}_\zeta(x,\overline{p}^L_\zeta(\varphi^i_{r,s}(F)))
=(-1)^{r+s}\zeta_i^{-(r(r-1)+s(s-1))/2}
\overline{\tau}^{\emptyset,L}_\zeta(\pi_\zeta(A_i)^sx\pi_\zeta(A_i)^r,\overline{p}^L_\zeta(F))
\]
for any $x\in U_\zeta^+$.
It follows that 
\begin{align*}
&\overline{\tau}^{\emptyset,L}_1(x',\xi^L(\overline{p}^L_\zeta(\varphi^i_{r,s}(F))))
=
\overline{\tau}^{\emptyset,L}_\zeta({}^t\xi(x'),\overline{p}^L_\zeta(\varphi^i_{r,s}(F)))\\
=&
(-1)^{r+s}\zeta_i^{-(r(r-1)+s(s-1))/2}
\overline{\tau}^{\emptyset,L}_\zeta(\pi_\zeta(A_i^s){}^t\xi(x')\pi_\zeta(A_i^r),\overline{p}^L_\zeta(F))\\
=&
(-1)^{r+s}\zeta_i^{-(r(r-1)+s(s-1))/2}
\overline{\tau}^{\emptyset,L}_\zeta({}^t\xi(x')\pi_\zeta(A_i^{r+s}),\overline{p}^L_\zeta(F))\\
=&\zeta_i^{rs}\overline{\tau}^{\emptyset,L}_1(x',\xi^L(\overline{p}^L_\zeta(\varphi^i_{r+s,0}(F))))
\end{align*}
for any $x'\in U^+_1$.
Here we have used the fact that the image of ${}^t\xi$ is contained in the center.
\end{proof}
\begin{lemma}
\label{lem:PU-}
Let $i\in I$, $\gamma\in Q^+$ and $b\in U^+_{1,\gamma}$.
\begin{itemize}
\item[\rm(i)]
Let $f\in \overline{U}^{L,-}_{1,-(\gamma+\alpha_i)}$.
Write
\[
\Delta(f)=\sum_{\gamma_1,\gamma_2\in Q^+,\gamma_1+\gamma_2=\gamma+\alpha_i}f_{\gamma_1,\gamma_2}\qquad
(f_{\gamma_1,\gamma_2}\in
\overline{U}^{L,-}_{1,-\gamma_1}\otimes\overline{U}^{L,-}_{1,-\gamma_2}),
\]
and define $f'\in \overline{U}^{L,-}_{1,-\gamma}$ by 
$f_{\gamma,\alpha_i}=f'\otimes f_i$.
Then we have
\begin{equation}
\label{lem:PU-1}
\overline{\tau}^{\emptyset,L}_1(\{\pi_1(A_i),b\},f)
=
\overline{\tau}_1^{\emptyset,L}
\left(b,
\frac{(\alpha_i,\alpha_i)}2[f,e_i]
-
\frac{(\alpha^\vee_i,\gamma)}2f'
\right).
\end{equation}
\item[\rm(ii)]
Let $f\in \overline{U}^{L,-}_{\zeta,-\ell(\gamma+\alpha_i)}$.
Write
\begin{multline*}
\Delta(f)=\sum_{\gamma_1,\gamma_2\in Q^+,\gamma_1+\gamma_2=\ell(\gamma+\alpha_i)}f_{\gamma_1,\gamma_2}
(1\otimes\overline{p}^L_\zeta(K_{-\gamma_1}))
\\
(f_{\gamma_1,\gamma_2}\in
\overline{U}^{L,-}_{\zeta,-\gamma_1}\otimes\overline{U}^{L,-}_{\zeta,-\gamma_2}),
\end{multline*}
and define $f'\in \overline{U}^{L,-}_{\zeta,-\ell\gamma}$ by 
$f_{\ell\gamma,\ell\alpha_i}=f'\otimes \overline{p}^L_\zeta(F_i^{(\ell)})$.
Then we have
\begin{multline}
\label{lem:PU-2}
\overline{\tau}^{\emptyset,L}_\zeta(\{\pi_\zeta(A_i^\ell),{}^t\xi(b)\},f)\\
=
\overline{\tau}_1^{\emptyset,L}
\left(b,
\frac{(\alpha_i,\alpha_i)}2[\xi^L(f),e_i]
-
\frac{(\alpha^\vee_i,\gamma)}2\xi^L(f')
\right).
\end{multline}
\end{itemize}
\end{lemma}
\begin{proof}
Note that (i)  can be regarded as a special case of (ii) when $\ell=1$.
Hence we will only prove (ii).
Take $F\in U^{L,-}_{\BA_\zeta,-\ell(\gamma+\alpha_i)}$ such that 
$\overline{p}^L_\zeta(F)=f$.
By
\[
\Delta_2(E_i^{(\ell)})
=\sum_{r+s+t=\ell}
q_i^{rs+st+tr}
E_i^{(r)}K_i^{s+t}\otimes
E_i^{(s)}K_i^t\otimes E_i^{(t)}
\]
and Lemma \ref{lem:Drinfeld paring} we have 
\begin{align*}
E_i^{(\ell)} F
=&\sum_{r+s+t=\ell}\sum_{(F)_2}
q_i^{rs+st+tr}
\tau(E_i^{(r)}K_i^{s+t},F_{(0)})\tau(E_i^{(t)},SF_{(2)})F_{(1)}E_i^{(s)}K_i^t\\
=&\sum_{r+s+t=\ell}\sum_{(F)_2}
q_i^{rs-st+tr}
\tau(E_i^{(r)},F_{(0)})\tau(E_i^{(t)},SF_{(2)})F_{(1)}K_i^tE_i^{(s)}\\
=&\sum_{s=0}^\ell X_sE_i^{(s)},
\end{align*}
where
\[
X_s=
\sum_{r+t=\ell-s}\sum_{(F)_2}
q_i^{rs-st+tr}
\tau(E_i^{(r)},F_{(0)})\tau(E_i^{(t)},SF_{(2)})F_{(1)}K_i^t.
\]
By $E_i^{(\ell)} F\in U_{\BA_\zeta}^L$ 
we have $X_s\in U_{\BA_\zeta}^L$ for $0\leqq s\leqq\ell$.
Note that $X_\ell=F$ and 
\[
X_0=\sum_{r+t=\ell}X_{0,r,t},
\]
where
\[
X_{0,r,t}=(-1)^t
q_i^{t(1-\ell-\ell(\alpha_i^\vee,\gamma))}
\tau(E_i^{(r)},F_i^{(r)})
\tau(E_i^{(t)},F_i^{(t)})
\varphi^i_{r,t}(F)K_i^{-r+t}.
\]
Note also that
\[
\tau(E_i^{(m)},F_i^{(m)})=
\frac{(-1)^mq_i^{m(m-1)/2}}
{[m]_{q_i}!(q_i-q_i^{-1})^m}.
\]
Hence $X_{0,r,t}\in U_{\BA_{\zeta}}^L$ for $r+t=\ell, \,r\ne0, \,t\ne0$.
By
$X_0\in U_{\BA_\zeta}^L$ we also have
$X_{0,\ell,0}+X_{0,0,\ell}\in  U_{\BA_{\zeta}}^L$.
From this we obtain
\[
\varphi^i_{\ell,0}(F)K_i^{-\ell}-
q_i^{\ell(1-\ell-\ell(\alpha_i^\vee,\gamma))}
\varphi^i_{0,\ell}(F)K_i^{\ell}
\in (q^\ell-q^{-\ell})U_{\BA_\zeta}^L,
\]
or equivalently,
\[
\varphi^i_{\ell,0}(F)K_i^{-\ell}-
\varphi^i_{0,\ell}(F)K_i^{\ell}
\in (q^\ell-q^{-\ell})U_{\BA_\zeta}^L.
\]
Let us show
\begin{equation}
\label{eq:hen}
\xi^L(\overline{p}^L_\zeta(\sum_{r+t=\ell, r>0,t>0}X_{0,r,t}))=
\frac{\ell-1}{2\ell}
\xi^L(\overline{p}^L_\zeta(\varphi_{\ell,0}^i(F)))
\end{equation}
By Lemma \ref{lem:varphi} we have
\begin{align*}
&\xi^L(\overline{p}^L_\zeta(\sum_{r+t=\ell, r>0,t>0}X_{0,r,t}))\\
=&
\sum_{t=1}^{\ell-1}
(-1)^t\zeta_i^t
\frac
{
(-1)^{\ell-t}\zeta_i^{(\ell-t)(\ell-t-1)/2}
}
{
[\ell-t]_{\zeta_i}!(\zeta_i-\zeta_i^{-1})^{\ell-t}
}
\frac
{
(-1)^{t}\zeta_i^{t(t-1)/2}
}
{
[t]_{\zeta_i}!(\zeta_i-\zeta_i^{-1})^{t}
}
\zeta_i^{t(\ell-t)}
\xi^L(\overline{p}^L_\zeta(\varphi_{\ell,0}^i(F)))\\
=&
-\frac1{(\zeta_i-\zeta_i^{-1})^{\ell}}
\sum_{t=1}^{\ell-1}
\frac
{
(-1)^t\zeta_i^t
}
{
[\ell-t]_{\zeta_i}![t]_{\zeta_i}!}\xi^L(\overline{p}^L_\zeta(\varphi_{\ell,0}^i(F)))
,
\end{align*}
and hence it is sufficient to show
\[
\sum_{t=1}^{\ell-1}
\frac
{
(-1)^t\zeta_i^t
}
{
[\ell-t]_{\zeta_i}![t]_{\zeta_i}!
}=(\zeta_i-\zeta_i^{-1})^\ell\frac{1-\ell}{2\ell}.
\]
Indeed we have
\begin{align*}
&2\sum_{t=1}^{\ell-1}
\frac{(-1)^t\zeta_i^t}{[\ell-t]_{\zeta_i}![t]_{\zeta_i}!}=
\sum_{t=1}^{\ell-1}
\frac{(-1)^t\zeta_i^t}{[\ell-t]_{\zeta_i}![t]_{\zeta_i}!}
+
\sum_{t=1}^{\ell-1}
\frac{(-1)^{\ell-t}\zeta_i^{\ell-t}}{[t]_{\zeta_i}![\ell-t]_{\zeta_i}!}\\
=&\sum_{t=1}^{\ell-1}(-1)^t
\frac{\zeta_i^t-\zeta_i^{-t}}{[\ell-t]_{\zeta_i}![t]_{\zeta_i}!}
=\frac{(\zeta_i-\zeta_i^{-1})}{[\ell-1]_{\zeta_i}!}
\sum_{t=1}^{\ell-1}(-1)^t
\begin{bmatrix}
\ell-1\\t-1
\end{bmatrix}_{\zeta_i}\\
=&\frac{(\zeta_i-\zeta_i^{-1})}{\ell}
\left(
1-\sum_{s=0}^{\ell-1}(-1)^s
\begin{bmatrix}
\ell-1\\s
\end{bmatrix}_{\zeta_i}
\right)
=
\frac{(\zeta_i-\zeta_i^{-1})}{\ell}
\left(1-
\prod_{j=0}^{\ell-2}
(\zeta_i^{\ell-2-2j}-1)
\right)\\
=&
\frac{(\zeta_i-\zeta_i^{-1})}{\ell}
\left(1-
\prod_{j=1}^{\ell-1}
(1-\zeta_i^{-2j})
\right)
=(\zeta_i-\zeta_i^{-1})^\ell
\frac{1-\ell}\ell.
\end{align*}
\eqref{eq:hen} is proved.
On the other hand we have
\begin{align*}
&\xi^L(\overline{p}^L_\zeta(X_{0,\ell,0}+X_{0,0,\ell}))\\
=&
\xi^L
\left(\overline{p}^L_\zeta
\left(
-q_i^{\ell(\ell-1)/2}
\frac
{\varphi^i_{\ell,0}(F)K_i^{-\ell}-
q_i^{\ell(1-\ell-\ell(\alpha_i^\vee,\gamma))}
\varphi^i_{0,\ell}(F)K_i^{\ell}}
{[\ell]_{q_i}!(q_i-q_i^{-1})^\ell}
\right)\right)\\
=&
-\xi^L
\left(\overline{p}^L_\zeta
\left(
\frac
{\varphi^i_{\ell,0}(F)K_i^{-\ell}-
\varphi^i_{0,\ell}(F)K_i^{\ell}}
{[\ell]_{q_i}!(q_i-q_i^{-1})^\ell}
\right)\right)\\
&\qquad\qquad-
\xi^L
\left(\overline{p}^L_\zeta
\left(
\frac
{
(1-q_i^{\ell(1-\ell-\ell(\alpha_i^\vee,\gamma))})
\varphi^i_{0,\ell}(F)K_i^{\ell}}
{[\ell]_{q_i}!(q_i-q_i^{-1})^\ell}
\right)\right)\\
=&
-\xi^L
\left(\overline{p}^L_\zeta
\left(
\frac
{
\varphi^i_{\ell,0}(F)K_i^{-\ell}-
\varphi^i_{0,\ell}(F)K_i^{\ell}
}
{
\ell(q_i^\ell-q_i^{-\ell})
}
\right)\right)\\
&\qquad\qquad
+
\frac{1-\ell-\ell(\alpha_i^\vee,\gamma)}{2\ell}
\xi^L
\left(\overline{p}^L_\zeta
\left(
\varphi^i_{0,\ell}(F)
\right)\right),
\end{align*}
and hence
\begin{align*}
&[\xi^L(f),e_i]=
\xi^L(\overline{p}^L_\zeta([F,E_i^{(\ell)}])=
-\sum_{t=0}^\ell
\xi^L
(\overline{p}^L_\zeta(X_{0,\ell-t,t}))
\\
=&\xi^L
\left(\overline{p}^L_\zeta
\left(
\frac
{
\varphi^i_{\ell,0}(F)K_i^{-\ell}-
\varphi^i_{0,\ell}(F)K_i^{\ell}
}
{
\ell(q_i^\ell-q_i^{-\ell})
}
\right)\right)
+
\frac{(\alpha_i^\vee,\gamma)}{2}
\xi^L
\left(\overline{p}^L_\zeta
\left(
\varphi^i_{0,\ell}(F)\right)\right).
\end{align*}

Write
\begin{multline*}
\Delta(F)=\sum_{\gamma_1,\gamma_2\in Q^+,\gamma_1+\gamma_2=\ell(\gamma+\alpha_i)}F_{\gamma_1,\gamma_2}
(1\otimes K_{-\gamma_1})
\\
(F_{\gamma_1,\gamma_2}\in
{U}^{L,-}_{\BA,-\gamma_1}\otimes{U}^{L,-}_{\BA,-\gamma_2}).
\end{multline*}
Then we have
\[
F_{\ell\alpha_i,\ell\gamma}=
F_i^{(\ell)}\otimes \varphi^i_{\ell,0}(F),\qquad
F_{\ell\gamma,\ell\alpha_i}=
\varphi^i_{0,\ell}(F)\otimes F_i^{(\ell)}.
\]
Take ${B}\in U^+_{\BA_\zeta, \ell\gamma}$ such that $\pi_\zeta({B})={}^t\xi(b)$.
Then we have
\begin{align*}
&\overline{\tau}^{\emptyset,L}_\zeta(\{\pi_\zeta(A_i^\ell),{}^t\xi(b)\},f)\\
=&(\tau(A_i^\ell {B}-{B}A_i^\ell,F)/\ell(q^\ell-q^{-\ell}))|_{q=\zeta}\\
=&((\tau\otimes\tau)(B\otimes A_i^\ell -A_i^\ell\otimes B,\Delta(F))/\ell(q^\ell-q^{-\ell}))|_{q=\zeta}\\
=&(-q_i^{\ell(\ell-1)/2}\tau(B,\varphi^i_{0,\ell}(F)-\varphi^i_{\ell,0}(F))/\ell(q^\ell-q^{-\ell}))
|_{q=\zeta}\\
=&(\tau(B,\varphi^i_{\ell,0}(F)-\varphi^i_{0,\ell}(F))/\ell(q^\ell-q^{-\ell}))
|_{q=\zeta}\\
=&
(\tau({B},\varphi^i_{\ell,0}(F)K_i^{-\ell}-\varphi^i_{0,\ell}(F)K_i^{\ell})/\ell(q^\ell-q^{-\ell}))
|_{q=\zeta}\\
=&
\overline{\tau}^{\emptyset,L}_\zeta
\left({}^t\xi(b),
\overline{p}^L_\zeta
\left(
\frac
{\varphi^i_{\ell,0}(F)K_i^{-\ell}-\varphi^i_{0,\ell}(F)K_i^{\ell}}{\ell(q^\ell-q^{-\ell})}
\right)\right)
\\
=&\frac{(\alpha_i,\alpha_i)}2
\overline{\tau}^{\emptyset,L}_1
\left(b,
\xi^L\left(
\overline{p}^L_\zeta
\left(
\frac
{\varphi^i_{\ell,0}(F)K_i^{-\ell}-\varphi^i_{0,\ell}(F)K_i^{\ell}}{\ell(q_i^\ell-q_i^{-\ell})}
\right)\right)\right)
\\
=&\frac{(\alpha_i,\alpha_i)}2
\overline{\tau}^{\emptyset,L}_1
\left(b,
[\xi^L(f),e_i]-\frac{(\alpha_i^\vee,\gamma)}2
\xi^L\left(
\overline{p}^L_\zeta
\left(
\varphi^i_{0,\ell}(F)
\right)\right)\right)
\\
=&\frac{(\alpha_i,\alpha_i)}2
\overline{\tau}^{\emptyset,L}_1
\left(b,
[\xi^L(f),e_i]-\frac{(\alpha_i^\vee,\gamma)}2
\xi^L\left(
f'
\right)\right).
\end{align*}
We are done.
\end{proof}

Now let us finish the proof of Theorem \ref{thm:Poisson}.
Regarding $\BC[M]$ as a subspace of $U(\Gm)^*$
the Poisson bracket of $\BC[M]$ is uniquely determined by the following properties.
\begin{align}
\label{eq:P1}
&\langle\{\varphi,\psi\},1\rangle=0
\qquad(\varphi,\psi\in\BC[M]),\\
\label{eq:P2}
&\langle\{\varphi,\psi\},a\rangle=
\tilde{\kappa}([(d\varphi)_e,(d\psi)_e],a)
\qquad(\varphi,\psi\in\BC[M], a\in\Gm),\\
\label{eq:P3}
&\langle\{\varphi,\psi\},uv\rangle\\
&\qquad=
\sum_{(\varphi), (\psi)}
\left\langle
\varphi_{(0)}\psi_{(0)}\otimes\{\varphi_{(1)},\psi_{(1)}\}+
\{\varphi_{(0)},\psi_{(0)}\}\otimes\varphi_{(1)}\psi_{(1)},
u\otimes v\right\rangle\nonumber\\
&\qquad\qquad
(\varphi,\psi\in\BC[M], u, v\in U(\Gm)).\nonumber
\end{align}
In \eqref{eq:P2} we have identified the cotangent space $\Gm^*$ of $M$ at the identity element $e$ with the Lie algebra $\Gk$ via $\tilde{\kappa}$.
Define $\res:\BC[M]\to\Gk$ as the composite of 
\[
\BC[M]\hookrightarrow U(\Gm)^*\to\Gm^*\simeq\Gk,
\]
where $U(\Gm)^*\to\Gm^*$ is induced by the canonical embedding $\Gm\hookrightarrow U(\Gm)$ and $\Gm^*\simeq\Gk$ is given by $\tilde{\kappa}$.
Then we have $(d\varphi)_e=\res(\varphi)$ for any $\varphi\in\BC[M]$.
Therefore, the Poisson bracket of $\BC[M]$ is uniquely characterized as the bilinear map
\begin{equation}
\label{eq:PH}
\{\,,\,\}:\BC[M]\times\BC[M]\to\BC[M]
\end{equation}
satisfying 
\begin{itemize}
\item[(a)]
$\BC[M]$ becomes a Poisson Hopf algebra,
\item[(b)]
$\res(\{\varphi,\psi\})=[\res(\varphi),\res(\psi)]
\quad(\varphi, \psi\in \BC[M]).$
\end{itemize}

Assume that we are given a bilinear map \eqref{eq:PH} satisfying the condition (a).
Set 
\[
D=\{a_i, b_i,\chi_\lambda\mid i\in I, \lambda\in Q\}
\subset\BC[M].
\]
By Lemma \ref{lem:P-gen} and the general formula
\begin{equation}
\label{eq:prod-Poisson}
\res(\varphi\psi)=\varepsilon(\varphi)\res(\psi)+\varepsilon(\psi)\res(\varphi)\qquad
(\varphi, \psi\in\BC[M])
\end{equation}
we see that the condition (b) is satisfied if and only if 
\begin{itemize}
\item[(c)]
$\res(\{\varphi,\psi\})=[\res(\varphi),\res(\psi)]
\quad(\varphi\in D, \psi\in \BC[M]).$
\end{itemize}
Note that $\BC[M]$ is generated as an algebra by $\BC[M^+]$, $\BC[M^-]$ and $\BC[M^0]$, where $\BC[M^0]$, $\BC[M^\pm]$ are regarded as subalgebras of $\BC[M]$ by \eqref{eq:M-docomp}.
Hence by \eqref{eq:prod-Poisson}
we see that the condition (c) is satisfied if and only if 
\begin{itemize}
\item[(d)]
$\res(\{\varphi,\psi\})=[\res(\varphi),\res(\psi)]
\quad(\varphi\in D, \psi\in \BC[M^+]\cup\BC[M^-]\cup\BC[M^0]).$
\end{itemize}

Let us show that the isomorphism $\Upsilon:U_1\to\BC[M]$ preserves the Poisson structure.
Set $r=\res\circ\Upsilon:U_1\to\Gk$, and
\begin{equation}
\label{Dprime}
D'=\{\pi_1({A}_i), \pi_1({B}_i), \pi_1({K}_\lambda)\mid i\in I, \lambda\in Q\}\subset U_1.
\end{equation}
By the above argument it is sufficient to show
\begin{equation}
\label{eq:Poisson}
r(\{a,b\})=[r(a),r(b)]
\qquad(a\in D', b\in U_1^+\cup U_1^-\cup U_1^0),
\end{equation}
where the Poisson bracket $\{\,,\,\}$ of $U_1$ is given by \eqref{eq:PoissonU1}.
We have
\begin{align*}
&r(\pi_1({A}_i))=\frac{(\alpha_i,\alpha_i)}2\theta(e_i),\qquad
r(\pi_1({B}_i))=-\frac{(\alpha_i,\alpha_i)}2\theta(f_i),\\
&r(\pi_1({K}_\lambda))=\frac12\theta(h_\lambda),
\end{align*}
where $h_\lambda\in\Gh$ is such that $\kappa(h_\lambda,h)=\lambda(h)$ for any $h\in\Gh$.

Let us show \eqref{eq:Poisson} for $a=\pi_1({K}_\lambda)$.
If $b=\pi_1({B})$ for $B\in U_\BA^0$, then we have $\{a,b\}=0$ by $[U_\BA^0,U_\BA^0]=0$, and $[r(a),r(b)]=0$ by $r(a), r(b)\in\theta(\Gh)$.
Assume $b=\pi_1({B})$ for $B\in U_{\BA,\pm\gamma}^\pm$.
Then we have
\begin{align*}
\{a,b\}
=&\pi_1((K_\lambda B-B K_\lambda)/(q-q^{-1}))
=\pi_1((q^{\pm(\lambda,\gamma)}-1)B K_\lambda/(q-q^{-1}))\\
=&\pm(\lambda,\gamma)b/2.
\end{align*}
On the other hand by $r(b)\in\theta\Gg_{\pm\gamma}$ we have 
$[r(a), r(b)]=\pm(\lambda,\gamma)r(b)/2$.
Hence \eqref{eq:Poisson} is proved for $a=\pi_1({K}_\lambda)$.
Note that the above argument also give the proof for the case $b\in U_1^0$.
It remains to show  \eqref{eq:Poisson} when $a=\pi_1({A}_i)$ or $\pi_1({B}_i)$, and $b\in U_1^\pm$.

Let us consider the case $a=\pi_1({A}_i)$ and $b\in U_{1,-\gamma}^-$ for $\gamma\in Q^+$.
If $\gamma-\alpha_i\notin\Delta\cup\{0\}$, then both sides of \eqref{eq:Poisson} are zero by $r(\{a,b\}), [r(a),r(b)]\in\theta(\Gg_{-\gamma+\alpha_i})$.
In the case $\gamma=0$ we can easily check that the both sides of \eqref{eq:Poisson} coincide with zero.
In the case $\gamma=\alpha_i$ we can also easily check that the both sides of \eqref{eq:Poisson} coincide with $(\alpha_i,\alpha_i)\theta(h_{\alpha_i})/2$.
Therefore, we may assume that $\gamma-\alpha_i\in \Delta^+$.
In this case it is sufficient to show
\[
\tilde{\kappa}(r(\{a,b\}),(x,0))=
\tilde{\kappa}([r(a),r(b)],(x,0))
\]
for any $x\in\Gg_{\gamma-\alpha_i}$.
By Lemma \ref{lem:PU+} we have
\[
\{a,b\}
=\frac{(\alpha_i,\alpha_i)}2(b''\pi_1(K_i)-b'\pi_1(K_i^{-1})),
\]
where $b', b''$ are as in \eqref{eq:PU+}, and hence
\begin{align*}
&\tilde{\kappa}(r(\{a,b\}),(x,0))\\
=&
\frac{(\alpha_i,\alpha_i)}2
\tilde{\kappa}(r(b''\pi_1(K_i)-b'\pi_1(K_i^{-1})),(x,0))\\
=&
\frac{(\alpha_i,\alpha_i)}2
\overline{\sigma}_1(b''\pi_1(K_i)-b'\pi_1(K_i^{-1}),(x,0))\\
=&
\frac{(\alpha_i,\alpha_i)}2
\overline{\sigma}_1(b''-b',(x,0)).
\end{align*}
On the other hand we have
\begin{align*}
&\tilde{\kappa}([r(a),r(b)],(x,0))\\
=&-\tilde{\kappa}(r(b),[r(a),(x,0)])\\
=&-(\alpha_i,\alpha_i)\tilde{\kappa}(r(b),[x_i,(x,0)])/2\\
=&-(\alpha_i,\alpha_i)\overline{\sigma}_1(b,[x_i,(x,0)])/2\\
=&-(\alpha_i,\alpha_i)
(\overline{\sigma}_1\otimes\overline{\sigma}_1)
(\Delta(b),x_i\otimes(x,0)-(x,0)\otimes x_i)/2\\
=&-(\alpha_i,\alpha_i)
(\overline{\sigma}_1(\pi_1(B_i),x_i)
\overline{\sigma}_1(b'\pi_1(K_i^{-1}),(x,0))
-
\overline{\sigma}_1(b'',(x,0))
\overline{\sigma}_1(\pi_1(B_iK_{-(\gamma-\alpha_i)}),x_i)
/2\\
=&-(\alpha_i,\alpha_i)
(\overline{\sigma}_1(b'-b'',(x,0))
/2.
\end{align*}
\eqref{eq:Poisson} is proved in the case $a=\pi_i(A_i)$ and $b\in U_1^-$.

Let us next consider the case $a=\pi_1({A}_i)$ and $b\in U^+_{1,\gamma}$.
We have $r(a)\in\theta(\Gg_{\alpha_i})$, $r(b)\in\theta(\Gg_{\gamma})$, $r(\{a,b\})\in\theta(\Gg_{\gamma+\alpha_i})$, and hence we may assume that $\gamma+\alpha_i\in\Delta^+$.
If $\gamma=0$, then the both sides of \eqref{eq:Poisson} is zero.
Hence we may also assume that $\gamma\in\Delta^+$.
Then it is sufficient to show 
\[
\tilde{\kappa}(r(\{a,b\}),(0,y))=\tilde{\kappa}([r(a), r(b)],(0,y))
\]
for $y\in\Gg_{-(\gamma+\alpha_i)}$.
By Lemma \ref{lem:PU-} we have
\[
\tilde{\kappa}(r(\{a,b\}),(0,y))
=
\overline{\sigma}_1(\{a,b\},(0,y))
=\overline{\tau}_1^{\emptyset,L}
(\{a,b\},y)
=\frac{(\alpha_i,\alpha_i)}2
\overline{\tau}_1^{\emptyset,L}
(b,[y,e_i]).
\]
On the other hand we have 
\begin{align*}
&\tilde{\kappa}([r(a), r(b)],(0,y))
=-\tilde{\kappa}(r(b),[r(a),(0,y)])\\
=&-(\alpha_i,\alpha_i)\tilde{\kappa}(r(b),(0,[e_i,y]))/2
=(\alpha_i,\alpha_i)
\overline{\sigma}_1(b,[y,e_i])/2\\
=&(\alpha_i,\alpha_i)
\overline{\tau}_1^{\emptyset,L}(b,[y,e_i])/2.
\end{align*}
\eqref{eq:Poisson} is proved in the case $a=\pi_i(A_i)$ and $b\in U_1^+$.
The remaining case $a=\pi_1(B_i)$ is proved similarly to the case $a=\pi_1(A_i)$. 
We omit the details.
Now we have proved that the isomorphism $\Upsilon:U_1\to\BC[M]$ preserves the Poisson structure.

It remains to show that $Z_\zeta$ is closed under the Poisson bracket and that the isomorphism ${}^t\xi:U_1\simeq Z_\zeta$ preserves the Poisson structure.
By the above argument and by Lemma \ref{lem:P-gen} we see that
the set $D'$ (see \eqref{Dprime}) generates the Poisson algebra.
Therefore, it is sufficient to show 
\begin{equation}
\label{eq:1zeta}
\{{}^t\xi(a), {}^t\xi(b)\}={}^t\xi(\{a,b\})\qquad
(a\in D', b\in U_1^+\cup U_1^0\cup U_1^-).
\end{equation}
The case $a=\pi_1(K_\lambda)$ is easy.
Hence it is sufficient to show \eqref{eq:1zeta} in the cases
$a=\pi_1(A_i)$ or $a=\pi_1(B_i)$, and $b\in U_1^+\cup U_1^-$.
Assume $a=\pi_1(A_i)$.
If $b\in U_1^-$, the assertion follows from Lemma \ref{lem:PU+}.
If  $b\in U_1^+$, then we see easily by Lemma \ref{lem:PU-} that
\[
\overline{\sigma}_\zeta(\{{}^t\xi(a),{}^t\xi(b)\},y)
=\overline{\sigma}_\zeta({}^t\xi(\{a,b\}),y)
\]
for any $y\in \overline{V}_\zeta^{-}$.
Since $\{{}^t\xi(a),{}^t\xi(b)\}, {}^t\xi(\{a,b\}\in U_\zeta^+$, 
\eqref{eq:1zeta} holds in this case.
The proof for the case $a=\pi_1(B_i)$ is similar to that for $a=\pi_1(A_i)$.
Details are omitted.
The proof of Theorem \ref{thm:Poisson}
is now complete.

%%%%%%%%%%%%%%%%%%%%%%%%%%%%%%%%
\section{Poisson manifold associated to rings of differential operators}
We denote by $F$ the subspace of $U^*$ spanned by the matrix coefficients of finite dimensional $U$-modules $E$ such that 
\[
E=\bigoplus_{\lambda\in Q}E_\lambda
\quad\mbox{ with }\quad
E_\lambda=\{v\in E\mid K_\mu v=q^{(\lambda,\mu)}v\,\,(\forall \mu\in Q)\}.
\]
It is endowed with a structure of Hopf algebra via
\begin{align*}
&\langle \varphi\psi,u\rangle=
\langle\varphi\otimes\psi,\Delta(u)\rangle
\qquad
&(\varphi,\psi\in F,\,\,u\in U),\\
&\langle 1,u\rangle=
\varepsilon(u)
\qquad
&(u\in U),\\
&\langle \Delta(\varphi),u\otimes u'\rangle=
\langle\varphi,uu'\rangle
\qquad
&(\varphi\in F,\,\,u, u'\in U),\\
&\epsilon(\varphi)=
\langle\varphi,1\rangle,
\qquad
&(\varphi\in F),\\
&\langle S(\varphi),u\rangle=
\langle\varphi,S(u)\rangle
\qquad
&(\varphi\in F,\,\,u\in U),
\end{align*}
where $\langle\,\,,\,\,\rangle:F\times U\to\BF$ is the canonical paring.
$F$ is also endowed with a structure of $U$-bimodule by
\[
\langle u'\varphi u'',u\rangle=\langle \varphi, u''uu'\rangle
\qquad
(\varphi\in F, u, u', u''\in U).
\]
For a subring $\BA$ of $\BF$ containing $\BQ[q,q^{-1}]$ we set
\[
F_{\BA}=
\{\varphi\in F\mid
\langle\varphi,U^L_{\BA}\rangle\subset\BA\}.
\]
It is a Hopf algebra over $\BA$ and a $U^L_{\BA}$-bimodule.
For $z\in\BC^\times$ we set
\[
F_z=\BC\otimes_{\BA_z}F_{\BA_z},
\]
where ${\BA_z}\to\BC$ is given by $q\mapsto z$.
Then $F_z$ is a Hopf algebra over $\BC$ and a $U^L_z$-bimodule.
In the following 
we will only be concerned with $F_1$, which is canonically isomorphic to the coordinate algebra $\BC[G]$ of the adjoint group of $\Gg$.

Denote by $\varphi\mapsto\overline{\varphi}$ the canonical homomorphism $F_{\BA_1}\to F_1=\BC[G]$.
We have a natural Poisson Hopf algebra structure on $\BC[G]$ given by
\[
\{\overline{\varphi},\overline{\psi}\}=\overline{[\varphi,\psi]/(q-q^{-1})}\qquad(\varphi,\psi\in F_{\BA_1}).\]
It is known that this Poisson Hopf algebra structure of $\BC[G]$ coincides with the one comming from the Manin triple $(\Gg\oplus\Gg,\Gm,\Gk)$ by identifying $\Gk$ with $\Gg$ (see De Concini-Lyubashenko) \cite{DL}).

We define a $\BF$-algebra structure on
\[
D=F\otimes_{\BF} U
\]
by
\[
(\varphi\otimes u)(\varphi'\otimes u')
=\sum_{(u)}\varphi(u_{(0)}\varphi')\otimes u_{(1)}u'\qquad
(\varphi, \varphi'\in F, u, u'\in U).
\]
The algebra $D$ is  an analogue of the ring of differential operators on $G$.
We will identify $U$ and $F$ with subalgebras of $D$ by the embeddings
$U\ni u\mapsto 1\otimes u\in D$ and $F\ni\varphi\mapsto\varphi\otimes1\in D$ respectively.

Let $\BA$ be a subring of $\BF$ containing $\BQ[q,q^{-1}]$.
We have a natural $\BA$-form
\[
D'_\BA=F_\BA\otimes_{\BA} U^L_\BA
\]
of $D$ whose specialization 
\[
D'_1=\BC\otimes_\BA D'_\BA=F_1\otimes U^L_1
\]
at $q=1$ is almost isomorphic to the ring $\BC[G]\otimes_\BC U(\Gg)
$ of differential operators on $G$.
However, in the following we will be concerned with a different $\BA$-form
\[
D_\BA=F_\BA\otimes_{\BA} U_\BA.
\]
For $z\in \BC^\times$ we set
\[
D_z=\BC\otimes_{\BA_z} D_{\BA_z}=F_z\otimes U_z
\]
where $\BA_z\to\BC$ is given by $q\mapsto z$.
\begin{lemma}
\label{lem:com-D}
$D_1$ is a commutative algebra.
In particular, it is identified as an algebra with the coordinate algebra $\BC[G]\otimes\BC[M]$ of $G\times M$.
\end{lemma}
\begin{proof}
By the definition of the multiplication of $D$ it is sufficient to show
$u\varphi=\varepsilon(u)\varphi$ for $u\in U_1$, $\varphi\in F_1$.
Let $\iota:U_1\to \overline{U}^L_1$ be the algebra homomorphism induced by $U_\BA\subset U^L_\BA$.
Then we have
\[
\langle u\varphi,u'\rangle
=\langle\varphi, u'\iota(u)\rangle 
\qquad(u'\in \overline{U}^L_1),
\]
and hence it is sufficient to show $\iota(u)=\varepsilon(u)1$ for any $u\in U_1$.
We may assume that $u$ is one of $K_\lambda\,\,(\lambda\in Q)$, $A_{\beta_k}, B_{\beta_k}\,\,(1\leqq k\leqq N)$.
In these cases the assertion follows from $\iota(K_\lambda)=1$, $\iota(A_{\beta_k})=\iota(B_{\beta_k})=0$
\end{proof}

\begin{remark}{\rm
We can show that $D_1$ is isomorphic to a central subalgebra of $D_\zeta$, where 
$\zeta$ is as in Section \ref{section:roots}.}
\end{remark}

By Lemma \ref{lem:com-D} we have a natural Poisson algebra structure of $D_1=\BC[G]\otimes\BC[M]$ given by
\[
\{\overline{\Phi},\overline{\Phi'}\}=
\overline{[\Phi,\Phi']/(q-q^{-1})}
\qquad(\Phi,\Phi'\in D_\BA),
\]
where $D_\BA\ni\Phi\mapsto\overline{\Phi}\in D_1=\BC[G]\otimes\BC[M]$ is the natural homomorphism.
Let us describe this Poisson bracket more explicitly.

By definition the canonical inclusions $\BC[G]\ni\varphi\mapsto\varphi\otimes 1\in D_1$ and $\BC[M]\ni\psi\mapsto1\otimes\psi\in D_1$ are homomorphisms of Poisson algebras.
Since the Poisson structures of $\BC[G]$ and $\BC[M]$ are already described explicitly, we have only to give a description of $\{\varphi,\psi\}$ for $\varphi\in\BC[G]$, $\psi\in\BC[M]$.

In general, for an algebraic group $S$ with Lie algebra $\Gs$ we denote by 
\[
\langle\,\,,\,\,\rangle:\BC[S]\times U(\Gs)\to\BC
\]
the canonical Hopf paring.
We have a $U(\Gs)$-bimodule structure of $\BC[S]$ given by 
\[
\langle u'\varphi u'',u\rangle=\langle\varphi,u''uu'\rangle.
\]
For $a\in\Gs, \varphi\in\BC[S], s\in S$ we have
\[
(a\varphi)(s)=\frac{d}{dt}\varphi(s\exp(ta))|_{t=0},\qquad
(\varphi a)(s)=\frac{d}{dt}\varphi(\exp(ta)s)|_{t=0}.
\]
For $a\in\Gs$ we denote by $L_a$ (resp. $R_a$) the left (resp. right) invariant vector field on $S$ given by $L_a(\varphi)=a\varphi$ (resp. $R_a(\varphi)=\varphi a$).
For $b\in\Gs^*$ we denote by $L^*_b$ (resp. $R^*_b$) the  left (resp. right) invariant 1-form on $S$ given by $\langle L_a,L^*_b\rangle=\langle a,b\rangle$ (resp. $\langle R_a,R^*_b\rangle=\langle a,b\rangle$).

\begin{proposition}
\label{prop:GMPoisson}
For $\varphi\in\BC[G], \psi\in\BC[M]$ we have
\[
\{\varphi,\psi\}
=-\sum_{r=1}^{\dim\Gg}(L_{\xi_r}(\varphi))(R_{\eta_r}(\psi)),
\]
where $\{\xi_r\}_{r=1}^{\dim\Gg}$ and $\{\eta_r\}_{r=1}^{\dim\Gg}$ are bases of $\Gg$ and $\Gm$ respectively such that
$\tilde{\kappa}(\theta(\xi_r),\eta_s)=\delta_{rs}$.
\end{proposition}
\begin{proof}
Our assertion is equivalent to the identity
\[
\overline{(\sum_{(x)}(x_{(0)}f)\otimes x_{(1)}-f\otimes x)/(q-q^{-1})}
=\sum_r\xi_r\overline{f}\otimes\overline{x}\eta_r
\quad(f\in F_{\BA_1}, x\in U_{\BA_1})
\]
in $D_1=F_1\otimes U_1\cong\BC[G]\otimes \BC[M]$.
For $u\in U^L_{\BA_1}$ we have
\begin{align*}
&\sum_{(x)}\langle x_{(0)}f,u\rangle x_{(1)}-\langle f,u\rangle x
=\sum_{(x)}\langle f,ux_{(0)}\rangle x_{(1)}-\langle f,u\rangle x,\\
&
\sum_r\langle\xi_r\overline{f},\overline u\rangle\overline{x}\eta_r
=\sum_r\langle\overline{f},\overline u\xi_r\rangle\overline{x}\eta_r,
\end{align*}
and hence our assertion is further equivalent to the identity
\[
\overline{(\sum_{(x)}ux_{(0)}\otimes x_{(1)}-u\otimes x)/(q-q^{-1})}
=\sum_r\overline{u}\xi_r\otimes\overline{x}\eta_r
\quad(u\in U^L_{\BA_1}, x\in U_{\BA_1})
\]
in $\overline{U}^L_1\otimes U_1\cong U(\Gg)\otimes\BC[M]$.
This statement follows from its special case $u=1$:
\[
\overline{(\sum_{(x)}x_{(0)}\otimes x_{(1)}-1\otimes x)/(q-q^{-1})}
=\sum_r\xi_r\otimes\overline{x}\eta_r
\quad(x\in U_{\BA_1}).
\]
Let $\{x_j\}_j$ be a free basis of $U_{\BA_1}$, and define $v_j\in U_{\BA_1}^*=\Hom_{\BA_1}(U_{\BA_1},\BA_1)$ by $\langle v_j,x_k\rangle=\delta_{jk}$.
Then for $v\in U_{\BA_1}^*$ we have
\begin{align*}
&
\sum_{(x)}\langle v,x_{(1)}\rangle x_{(0)}-\langle v,x\rangle1
=\sum_{(x),j}\langle v,x_{(1)}\rangle 
\langle v_j,x_{(0)}\rangle x_j-\langle v,x\rangle1
\\
=&\sum_{j}
\langle v_jv,x\rangle 
x_j-\langle v,x\rangle1
,\\
&
\sum_r\langle \overline{v},\overline{x}\eta_r\rangle\xi_r
=\sum_r\langle \eta_r\overline{v},\overline{x}\rangle\xi_r.
\end{align*}
Here, the multiplication of $U_{\BA_1}^*$ is induced by the comultiplication of $U_{\BA_1}$.
Therefore, we have only to show the identity 
\[
\overline{(\sum_{j}x_{j}\otimes v_jv-1\otimes v)/(q-q^{-1})}
=\sum_r\xi_r\otimes\eta_r\overline{v}
\quad(v\in U^*_{\BA_1})
\]
in the completion $\Hom_\BC(\BC[M],U(\Gg))$ of $U(\Gg)\otimes\BC[M]^*$.
This statement follwos from its special case $v=1_{U_{\BA_1}^*}=\varepsilon$:
\[
\overline{(\sum_{j}x_{j}\otimes v_j-1\otimes 1)/(q-q^{-1})}
=\sum_r\xi_r\otimes\eta_r.
\]
This follows from Lemma \ref{lem:final} below.
\end{proof}
\begin{lemma}
\label{lem:final}
Let $\Xi:\BC[M]\to U(\Gg)$ be the map induced by
\[
U_{\BA_1}\ni x\mapsto(x-\varepsilon(x)1)/(q-q^{-1})\in U^L_{\BA_1}.
\]
Then we have $\Image(\Xi)\subset\Gg$ and 
\[
\tilde{\kappa}(\theta(\Xi(\varphi)),\eta)=\langle\varphi,\eta\rangle
\qquad (\varphi\in\BC[M], \,\,\eta\in\Gm).
\]
\end{lemma}
This can be shown by a direct computation in terms of root vectors.
Details are omitted.
\begin{remark}
{\rm
In terms of the Poisson tensor $\delta$ of the Poisson manifold $G\times M$ Proposition \ref{prop:GMPoisson} can be reformulated as follows.
Under the identification $\Gg\cong\Gm^*$, $\Gm^*\cong\Gg$ via $\tilde{\kappa}$ we have
\[
\delta_{(g,m)}((L^*_\eta)_g,(R^*_\xi)_m)=-\tilde{\kappa}(\xi,\eta)
\qquad((g, m)\in G\times M, \eta\in\Gm, \xi\in\Gg).
\]
}
\end{remark}

%%%%%%%%%%%%%%%%%

\bibliographystyle{unsrt}

\end{document}